\documentclass[english,12pt]{amsart}
\usepackage{nicefrac}
\usepackage[inline]{enumitem}
\usepackage{enumitem}
\usepackage{amssymb}
\usepackage[usenames,dvipsnames,svgnames,table]{xcolor}
\definecolor{myRed}{rgb}{0.9,0.,.2}
\definecolor{myBlue}{rgb}{0.,0.,.6}
\definecolor{myGreen}{rgb}{0.1,0.7,0.1}
\definecolor{myViolet}{rgb}{102,0,153}
\usepackage[colorlinks=true,urlcolor=myBlue,linkcolor=myBlue,citecolor=myGreen,pagebackref=true]{hyperref}

\usepackage{fouriernc}
\usepackage[scaled=0.775]{helvet}
\usepackage{courier}
\usepackage[marginparwidth=2cm]{geometry}
\usepackage{array,float}

\usepackage{tikz}
\usetikzlibrary{shapes,arrows}
\usetikzlibrary{fit,positioning}
\usetikzlibrary{patterns,decorations.pathreplacing}
\usepackage{xkeyval}
\usepackage{moreverb}
\usepackage{epic}
\usepgfmodule{shapes,plot,decorations}

\usepackage{booktabs} 

\usepackage[normalem]{ulem} 

\newcommand{\N}{\mathbb{N}}
\newcommand{\Z}{\mathbb{Z}}

\newcommand{\R}{\mathbb{R}}

\newcommand{\C}{\mathcal{C}} 

\newcommand{\HH}{\mathbb{H}}

\newcommand{\Hb}{\mathbb{H}}

\newcommand{\PGL}{\mathrm{PGL}}

\newcommand{\G}{\Gamma}
\newcommand{\g}{\gamma}

\renewcommand{\C}{\mathcal{C}}

\renewcommand{\O}{\Omega}

\renewcommand{\S}{\mathbb{S}}

\newcommand{\Aut}{\textrm{Aut}}

\newcommand{\Isom}{\mathrm{Isom}}

\newcommand{\Hit}{\mathrm{Hit}}

\newcommand{\PS}{\mathbb{S}}
\newcommand{\SL}{\mathrm{SL}}
\newcommand{\GL}{\mathrm{GL}}

\usepackage{mathrsfs}  
\newcommand{\Cart}{\mathscr{A}}

\theoremstyle{plain}
\newtheorem{theorem}{Theorem}[section]

\newtheorem{proposition}[theorem]{Proposition}

\newtheorem{conj}[theorem]{Conjecture}
\newtheorem{ques}[theorem]{Question}

\theoremstyle{definition}

\newtheorem{definition}[theorem]{Definition}

\theoremstyle{remark}
\newtheorem{remark}[theorem]{Remark}

\title{Discrete Coxeter groups}

\author{Gye-Seon Lee}
\address{Department of Mathematics, Sungkyunkwan University, 2066 Seobu-ro, Jangan-gu, Suwon, 16419, South Korea}
\email{gyeseonlee@skku.edu}

\author{Ludovic Marquis}
\address{Univ Rennes, CNRS, IRMAR - UMR 6625, F-35000 Rennes, France}
\email{ludovic.marquis@univ-rennes1.fr}

\begin{document}
\setcounter{tocdepth}{1}

%
%
%



\begin{abstract}
Coxeter groups are a special class of groups generated by involutions. They play important roles in the various areas of mathematics. This survey particularly focuses on how one uses Coxeter groups to construct interesting examples of discrete subgroups of Lie groups.
\end{abstract}



\maketitle


\tableofcontents

\section{Introduction}

It is a fundamental problem in geometry and topology to understand discrete subgroups of Lie groups $G$. For example, when $G$ is the isometry group $\mathrm{Isom}(\Hb^d)$ of the hyperbolic $d$-space $\Hb^d$, the study of discrete subgroups of $\mathrm{Isom}(\Hb^d)$ is closely related to that of complete hyperbolic $d$-manifolds. More precisely, there is a one-to-one correspondence between torsion-free discrete subgroups $\Gamma$ of $\mathrm{Isom}(\Hb^d)$ and complete hyperbolic $d$-manifolds $\Hb^d/\Gamma$. 

\medskip

Convex cocompact subgroups of rank-one Lie groups $G$ are specially an important class of discrete subgroups of $G$. In particular, given a finitely generated group $\Gamma$, the space of representations $\rho : \Gamma \rightarrow G$ whose image is convex cocompact is \emph{open} in the representation space $\mathrm{Hom}(\Gamma,G)$, i.e., the space of all representations $\rho : \Gamma \rightarrow G$. So, if $\rho$ is convex cocompact and is not isolated, then all the nearby representations of $\rho$ are again convex cocompact, and hence discrete. 

\medskip

Recently, new notions of representations were introduced to generalize convex cocompact subgroups of rank-one Lie groups: Anosov representations in real semisimple Lie groups (see \cite{labourie_anosov,Guichard_Wienhard_Anosov}) and convex cocompact subgroups in real projective spaces (see \cite{anosov_on_RPV}). Such representations also have the property of openness. As new theories are developed, it is also important to have many examples to support them. From this perspective, the role of Coxeter groups is crucial. 

\medskip

The aim of this survey is to illustrate how one can build interesting examples of discrete Coxeter groups.

\subsection*{Acknowledgements}

G.-S. Lee was supported by the National Research Foundation of Korea(NRF) grant funded by the Korea government(MSIT) (No. 2020R1C1C1A01013667). L. Marquis acknowledges support by the Centre Henri Lebesgue (ANR-11-LABX-0020 LEBESGUE). 

We would like to thank Athanase Papadopoulos for carefully reading this paper and suggesting several improvements.

\section{Coxeter groups}

\subsection{What is a Coxeter group?}

A \emph{Coxeter matrix}\index{Coxeter matrix} $M=(m_{s,t})_{s,t \in S}$ on a finite set $S$ is a symmetric matrix with entries $m_{s,t} \in \{1,2, \dotsc, m, \dotsc,\infty \}$ such that the diagonal entries $m_{s,s}=1$ and off-diagonal entries $m_{s,t} \neq 1$. From any Coxeter matrix $M=(m_{s,t})_{s,t \in S}$, one may obtain the \emph{Coxeter group}\index{Coxeter group} $W$ of $M$ given by generators and relations:
$$
W = \big\langle s \in S \mid (s t)^{m_{s,t}}=1,\,\, \forall s,t \in S \big\rangle.
$$
Here, $(s t)^\infty = 1$ means that there is no relation between $s$ and $t$. Since $m_{s,s}=1$, each generator $s$ is an involution, i.e., $s^2 = 1$. We shall use the notation $W$, $W_S$ or $W_{S,M}$ for a Coxeter group, depending on what is important to stress. One should remember that a Coxeter group is a group with a preferred generating set, namely $S$. The \emph{rank}\index{Coxeter group!rank} of $W_S$ is the cardinality $\# S$ of $S$. 

\medskip

The \emph{Coxeter diagram}\index{Coxeter diagram} of the Coxeter group $W_S$ is the labeled graph $\mathscr{G}_{W}$ such that:
\begin{enumerate}[label=(\roman*)]
  \item the set of nodes\footnote{A Coxeter group often comes with a Coxeter polytope in such a way that the nodes of the Coxeter diagram are in bijection with the facets of the Coxeter polytope. We shall use the word \textit{node} of the Coxeter diagram rather than \textit{vertex} to make a distinction between the vertices of the Coxeter polytope and the vertices of the Coxeter diagram.} of $\mathscr{G}_{W}$ is the set $S$;
  \item two nodes $s,t\in S$ are connected by an edge $\overline{s t}$ of $\mathscr{G}_{W}$ if $m_{s,t} \in \{ 3, 4, \dotsc, \infty \}$;
  \item the label of the edge $\overline{s t}$ is $m_{s,t}$ if $m_{s,t} \in \{ 4, 5, \dotsc, \infty \}$.
\end{enumerate}

It is well-known that for any subset $T$ of $S$, the subgroup of $W_S$ generated by $T$ is the Coxeter group $W_{T,M'}$ with generating set $T$ and exponents $m'_{s,t} = m_{s,t}$ for every $s,t \in T$ (see \cite[Chap.\,IV, Th.~2]{Bourbaki_group_456}). Such a subgroup $W_{T}$ is called a \emph{standard subgroup}\index{Coxeter group!standard subgroup} of $W_S$. 

\medskip

The connected components of the graph $\mathscr{G}_{W_S}$ are graphs of the form $\mathscr{G}_{W_{S_i}}$, $i\in I$, where the $(S_i)_{i \in I}$ form a partition of $S$. The standard subgroups $W_{S_i}$ are called the \emph{irreducible components}\index{Coxeter group!irreducible component} of $W_S$. Since $m_{s,t} = 2$ if and only if $st=ts$, we see that the group $W_S$ is the direct product of the subgroups $W_{S_i}$ for $i\in I$. A Coxeter group $W_S$ is \emph{irreducible}\index{Coxeter group!irreducible} when the Coxeter diagram $\mathscr{G}_{W_S}$ is connected, i.e., $\# I = 1$. A subset $T$ of $S$ is said to be "\emph{something}" if the Coxeter group $W_T$ is "something". For example, the word "something" can be replaced by "irreducible", and so on. Two subsets $T, U \subset S$ are \emph{orthogonal}\index{orthogonal} if $m_{t,u}=2$ for every $t \in T$ and every $u \in U$. 

\medskip

The \emph{Cosine matrix}\index{Cosine matrix} of $W_{S,M}$ is the $S \times S$ symmetric matrix $\mathsf{C}_{W}$ whose entries are:
 $$\big( \mathsf{C}_{W} \big)_{s,t} = -2\cos \left( \frac{\pi}{m_{s,t}} \right) \quad \,\,\textrm{for every }  s,t \in S$$ 

\medskip

An irreducible Coxeter group $W$ is said to be \emph{spherical}\index{Coxeter group!spherical} (resp. \emph{affine}\index{Coxeter group!affine}) when the Cosine matrix $\mathsf{C}_{W}$ is positive definite (resp. positive semi-definite but not definite).

\begin{theorem}[{Coxeter \cite{coxeter_classi_simplex,classi_cox_poly_coxeter} and Margulis--Vinberg \cite{margu_vin}}]
Let $W_S$ be an irreducible Coxeter group. Then exactly one of the following is true:
\begin{enumerate}[label=(\roman*)]
	\item If $W_S$ is spherical, then $W_S$ is a finite group.
	\item If $W_S$ is affine, then $\# S \geqslant 2$ and $W_S$ is virtually\footnote{A group $G$ is \emph{virtually "something"}\index{virtually} if there is a finite index subgroup $H \leqslant G$ such that $H$ is "something".} $\Z^{\# S -1}$.
	\item Otherwise, $W_S$ is \emph{large}\index{large}, i.e., there exists a surjective homomorphism of a finite index subgroup of $W_S$ onto a free group on two generators.
\end{enumerate}
\end{theorem}

\begin{remark}
These three cases are clearly exclusive. Consequently, if an irreducible Coxeter group $W_S$ is finite (resp. infinite and virtually abelian), then $W_S$ is spherical (resp. affine).
\end{remark}

\begin{remark}
The irreducible spherical and irreducible affine Coxeter groups were classified by Coxeter \cite{coxeter_classi_simplex,classi_cox_poly_coxeter}; see also Witt \cite{classi_diag_witt}. The complete list can be found in Table \ref{tab:classi_spherical_affine}.
\end{remark}

A Coxeter group (not necessarily irreducible) is \emph{spherical}\index{Coxeter group!spherical} (resp. \emph{affine}\index{Coxeter group!affine}) when all its irreducible components are spherical (resp. affine).

\newcommand{\scalevalue}{0.43}

\begin{table}[h]
\centering
\begin{tabular}{m{1.9cm} m{3cm} m{0.7cm} m{1.7cm} m{3cm}}
{\footnotesize $I_2(p)$ $(p\geqslant 5)$}
&
\begin{tikzpicture}[thick,scale=\scalevalue, every node/.style={transform shape}]
\node[draw,circle] (I1) at (0,0) {};
\node[draw,circle,right=.8cm of I1] (I2) {};
\draw (I1) -- (I2)  node[above,midway] {\large $p$};
\end{tikzpicture}
&&
{\footnotesize $\widetilde{A}_1$}
&
\begin{tikzpicture}[thick,scale=\scalevalue, every node/.style={transform shape}]
\node[draw,circle] (I1) at (0,0) {};
\node[draw,circle,right=.8cm of I1] (I2) {};
\draw (I1) -- (I2)  node[above,midway] {\large $\infty$};
\end{tikzpicture}
\\
{\footnotesize $A_n$ $(n \geqslant 1)$}
&
\begin{tikzpicture}[thick,scale=\scalevalue, every node/.style={transform shape}]

\node[draw,circle] (A1) at (0,0) {};
\node[draw,circle,right=.8cm of A1] (A2) {};
\node[draw,circle,right=.8cm of A2] (A3) {};
\node[draw,circle,right=1cm of A3] (A4) {};
\node[draw,circle,right=.8cm of A4] (A5) {};

\draw (A1) -- (A2)  node[above,midway] {};
\draw (A2) -- (A3)  node[above,midway] {};
\draw[ dotted,thick] (A3) -- (A4) node[] {};
\draw (A4) -- (A5) node[above,midway] {};
\end{tikzpicture}
&&
{\footnotesize $\widetilde{A}_n$ $(n \geqslant 2)$}
&
\begin{tikzpicture}[thick,scale=\scalevalue, every node/.style={transform shape}]
\node[draw,circle] (A1) at (0,0) {};
\node[draw,circle,above=.8cm of A1] (A2) {};
\node[draw,circle,right=.8cm of A2] (A3) {};
\node[draw,circle,right=.8cm of A3] (A4) {};
\node[draw,circle,right=.8cm of A4] (A5) {};
\node[draw,circle,below=.8cm of A5] (A6) {};
\node[draw,circle,below=.8cm of A6] (A7) {};
\node[draw,circle,left=.8cm of A7] (A8) {};
\node[draw,circle,left=.8cm of A8] (A9) {};
\node[draw,circle,left=.8cm of A9] (A10) {};

\draw (A1) -- (A2)  node[above,midway] {};
\draw (A2) -- (A3)  node[above,midway] {};
\draw (A3) -- (A4) node[] {};
\draw (A4) -- (A5) node[above,midway] {};
\draw (A5) -- (A6) node[] {};
\draw (A6) -- (A7) node[] {};
\draw (A7) -- (A8) node[] {};
\draw[ dotted,thick] (A8) -- (A9) node[] {};
\draw (A9) -- (A10) node[] {};
\draw (A10) -- (A1) node[] {};
\end{tikzpicture}
\\
\\
{\footnotesize $B_n$ $(n \geqslant 2)$}
&
\begin{tikzpicture}[thick,scale=\scalevalue, every node/.style={transform shape}]
\node[draw,circle] (B1) at (0,0) {};
\node[draw,circle,right=.8cm of B1] (B2) {};
\node[draw,circle,right=.8cm of B2] (B3) {};
\node[draw,circle,right=1cm of B3] (B4) {};
\node[draw,circle,right=.8cm of B4] (B5) {};

\draw (B1) -- (B2)  node[above,midway] {\large $4$};
\draw (B2) -- (B3)  node[above,midway] {};
\draw[ dotted,thick] (B3) -- (B4) node[] {};
\draw (B4) -- (B5) node[above,midway] {};
\end{tikzpicture}
&&
{\footnotesize $\widetilde{B}_2=\widetilde{C}_2$}
&
\begin{tikzpicture}[thick,scale=\scalevalue, every node/.style={transform shape}]
\node[draw,circle] (H1) at (0,0) {};
\node[draw,circle,right=.8cm of H1] (H2) {};
\node[draw,circle,right=.8cm of H2] (H3) {};
\draw (H1) -- (H2)  node[above,midway] {\large $4$};
\draw (H2) -- (H3)  node[above,midway] {\large $4$};
\end{tikzpicture}
\\
{\footnotesize $H_3$}
&
\begin{tikzpicture}[thick,scale=\scalevalue, every node/.style={transform shape}]
\node[draw,circle] (H1) at (0,0) {};
\node[draw,circle,right=.8cm of H1] (H2) {};
\node[draw,circle,right=.8cm of H2] (H3) {};

\draw (H1) -- (H2)  node[above,midway] {\large $5$};
\draw (H2) -- (H3)  node[above,midway] {};
\end{tikzpicture}
&&
{\footnotesize $\widetilde{B}_n$ $(n \geqslant 3)$}
&
\begin{tikzpicture}[thick,scale=\scalevalue, every node/.style={transform shape}]
\node[draw,circle] (B1) at (0,0) {};
\node[draw,circle,right=.8cm of B1] (B2) {};
\node[draw,circle,right=.8cm of B2] (B3) {};
\node[draw,circle,right=1cm of B3] (B4) {};
\node[draw,circle,right=.8cm of B4] (B5) {};
\node[draw,circle,above right=.8cm of B5] (B6) {};
\node[draw,circle,below right=.8cm of B5] (B7) {};

\draw (B1) -- (B2)  node[above,midway] {\large $4$};
\draw (B2) -- (B3)  node[above,midway] {};
\draw[ dotted,thick] (B3) -- (B4) node[] {};
\draw (B4) -- (B5) node[above,midway] {};
\draw (B5) -- (B6) node[above,midway] {};
\draw (B5) -- (B7) node[above,midway] {};
\end{tikzpicture}
\\
{\footnotesize $H_4$}
&
\begin{tikzpicture}[thick,scale=\scalevalue, every node/.style={transform shape}]
\node[draw,circle] (HH1) at (0,0) {};
\node[draw,circle,right=.8cm of HH1] (HH2) {};
\node[draw,circle,right=.8cm of HH2] (HH3) {};
\node[draw,circle,right=.8cm of HH3] (HH4) {};

\draw (HH1) -- (HH2)  node[above,midway] {\large $5$};
\draw (HH2) -- (HH3)  node[above,midway] {};
\draw (HH3) -- (HH4)  node[above,midway] {};
\end{tikzpicture}
&&
{\footnotesize $\widetilde{C}_n$ $(n\geqslant 3)$}
&
\begin{tikzpicture}[thick,scale=\scalevalue, every node/.style={transform shape}]
\node[draw,circle] (C1) at (0,0) {};
\node[draw,circle,right=.8cm of C1] (C2) {};
\node[draw,circle,right=.8cm of C2] (C3) {};
\node[draw,circle,right=1cm of C3] (C4) {};
\node[draw,circle,right=.8cm of C4] (C5) {};

\draw (C1) -- (C2)  node[above,midway] {\large $4$};
\draw (C2) -- (C3)  node[above,midway] {};
\draw[ dotted,thick] (C3) -- (C4) node[] {};
\draw (C4) -- (C5) node[above,midway] {\large $4$};
\end{tikzpicture}
\\
\\
{\footnotesize $D_n$ $(n\geqslant 4)$}
&
\begin{tikzpicture}[thick,scale=\scalevalue, every node/.style={transform shape}]

\node[draw,circle] (D1) at (0,0) {};
\node[draw,circle,right=.8cm of D1] (D2) {};
\node[draw,circle,right=1cm of D2] (D3) {};
\node[draw,circle,right=.8cm of D3] (D4) {};
\node[draw,circle, above right=.8cm of D4] (D5) {};
\node[draw,circle,below right=.8cm of D4] (D6) {};

\draw (D1) -- (D2)  node[above,midway] {};
\draw[ dotted] (D2) -- (D3);
\draw (D3) -- (D4) node[above,midway] {};
\draw (D4) -- (D5) node[above,midway] {};
\draw (D4) -- (D6) node[below,midway] {};
\end{tikzpicture}
&&
{\footnotesize $\widetilde{D}_n$ $(n\geqslant 4)$}
&
\begin{tikzpicture}[thick,scale=\scalevalue, every node/.style={transform shape}]
\node[draw,circle] (D1) at (0,0) {};
\node[draw,circle,below right=0.8cm of D1] (D3) {};
\node[draw,circle,below left=0.8cm of D3] (D2) {};
\node[draw,circle,right=.8cm of D3] (DA) {};
\node[draw,circle,right=1cm of DA] (DB) {};
\node[draw,circle,right=.8cm of DB] (D4) {};
\node[draw,circle, above right=.8cm of D4] (D5) {};
\node[draw,circle,below right=.8cm of D4] (D6) {};

\draw (D1) -- (D3)  node[above,midway] {};
\draw (D2) -- (D3)  node[above,midway] {};
\draw (D3) -- (DA) node[above,midway] {};
\draw[ dotted] (DA) -- (DB);
\draw (D4) -- (DB) node[above,midway] {};
\draw (D4) -- (D5) node[above,midway] {};
\draw (D4) -- (D6) node[below,midway] {};
\end{tikzpicture}
\\
\\
{\footnotesize $F_4$}
&
\begin{tikzpicture}[thick,scale=\scalevalue, every node/.style={transform shape}]
\node[draw,circle] (F1) at (0,0) {};
\node[draw,circle,right=.8cm of F1] (F2) {};
\node[draw,circle,right=.8cm of F2] (F3) {};
\node[draw,circle,right=.8cm of F3] (F4) {};

\draw (F1) -- (F2)  node[above,midway] {};
\draw (F2) -- (F3)  node[above,midway] {\large $4$};
\draw (F3) -- (F4)  node[above,midway] {};
\end{tikzpicture}
&&
{\footnotesize $\widetilde{F}_4$}
&
\begin{tikzpicture}[thick,scale=\scalevalue, every node/.style={transform shape}]
\node[draw,circle] (F1) at (0,0) {};
\node[draw,circle,right=.8cm of F1] (F2) {};
\node[draw,circle,right=.8cm of F2] (F3) {};
\node[draw,circle,right=.8cm of F3] (F4) {};
\node[draw,circle,right=.8cm of F4] (F5) {};

\draw (F1) -- (F2)  node[above,midway] {};
\draw (F2) -- (F3)  node[above,midway] {\large $4$};
\draw (F3) -- (F4)  node[above,midway] {};
\draw (F4) -- (F5)  node[above,midway] {};
\end{tikzpicture}
\\
\\
&
&&
{\footnotesize $\widetilde{G}_2$}
&
\begin{tikzpicture}[thick,scale=\scalevalue, every node/.style={transform shape}]
\node[draw,circle] (HH1) at (0,0) {};
\node[draw,circle,right=.8cm of HH1] (HH2) {};
\node[draw,circle,right=.8cm of HH2] (HH3) {};

\draw (HH1) -- (HH2)  node[above,midway] {\large $6$};
\draw (HH2) -- (HH3)  node[above,midway] {};
\end{tikzpicture}
\\
\\
{\footnotesize $E_6$}

&
\begin{tikzpicture}[thick,scale=\scalevalue, every node/.style={transform shape}]
\node[draw,circle] (E1) at (0,0) {};
\node[draw,circle,right=.8cm of E1] (E2) {};
\node[draw,circle,right=.8cm of E2] (E3) {};
\node[draw,circle,right=.8cm of E3] (E4) {};
\node[draw,circle,right=.8cm of E4] (E5) {};
\node[draw,circle,below=.8cm of E3] (EA) {};

\draw (E1) -- (E2)  node[above,midway] {};
\draw (E2) -- (E3)  node[above,midway] {};
\draw (E3) -- (E4)  node[above,midway] {};
\draw (E4) -- (E5)  node[above,midway] {};
\draw (E3) -- (EA)  node[left,midway] {};
\end{tikzpicture}
&&
{\footnotesize $\widetilde{E}_6$}

&
\begin{tikzpicture}[thick,scale=\scalevalue, every node/.style={transform shape}]
\node[draw,circle] (E1) at (0,0) {};
\node[draw,circle,right=.8cm of E1] (E2) {};
\node[draw,circle,right=.8cm of E2] (E3) {};
\node[draw,circle,right=.8cm of E3] (E4) {};
\node[draw,circle,right=.8cm of E4] (E5) {};
\node[draw,circle,below=.8cm of E3] (EA) {};
\node[draw,circle,below=.8cm of EA] (EB) {};

\draw (E1) -- (E2)  node[above,midway] {};
\draw (E2) -- (E3)  node[above,midway] {};
\draw (E3) -- (E4)  node[above,midway] {};
\draw (E4) -- (E5)  node[above,midway] {};
\draw (E3) -- (EA)  node[left,midway] {};
\draw (EB) -- (EA)  node[left,midway] {};
\end{tikzpicture}
\\
{\footnotesize $E_7$}

&
\begin{tikzpicture}[thick,scale=\scalevalue, every node/.style={transform shape}]
\node[draw,circle] (EE1) at (0,0) {};
\node[draw,circle,right=.8cm of EE1] (EE2) {};
\node[draw,circle,right=.8cm of EE2] (EE3) {};
\node[draw,circle,right=.8cm of EE3] (EE4) {};
\node[draw,circle,right=.8cm of EE4] (EE5) {};
\node[draw,circle,right=.8cm of EE5] (EE6) {};
\node[draw,circle,below=.8cm of EE3] (EEA) {};

\draw (EE1) -- (EE2)  node[above,midway] {};
\draw (EE2) -- (EE3)  node[above,midway] {};
\draw (EE3) -- (EE4)  node[above,midway] {};
\draw (EE4) -- (EE5)  node[above,midway] {};
\draw (EE5) -- (EE6)  node[above,midway] {};
\draw (EE3) -- (EEA)  node[left,midway] {};
\end{tikzpicture}
&&
{\footnotesize $\widetilde{E}_7$}

&
\begin{tikzpicture}[thick,scale=\scalevalue, every node/.style={transform shape}]
\node[draw,circle] (EE1) at (0,0) {};
\node[draw,circle,right=.8cm of EE1] (EEB) {};
\node[draw,circle,right=.8cm of EEB] (EE2) {};
\node[draw,circle,right=.8cm of EE2] (EE3) {};
\node[draw,circle,right=.8cm of EE3] (EE4) {};
\node[draw,circle,right=.8cm of EE4] (EE5) {};
\node[draw,circle,right=.8cm of EE5] (EE6) {};
\node[draw,circle,below=.8cm of EE3] (EEA) {};

\draw (EE1) -- (EEB)  node[] {};
\draw (EE2) -- (EEB)  node[] {};
\draw (EE2) -- (EE3)  node[] {};
\draw (EE3) -- (EE4)  node[] {};
\draw (EE4) -- (EE5)  node[] {};
\draw (EE5) -- (EE6)  node[] {};
\draw (EE3) -- (EEA)  node[] {};
\end{tikzpicture}
\\
{\footnotesize $E_8$}

&
\begin{tikzpicture}[thick,scale=\scalevalue, every node/.style={transform shape}]
\node[draw,circle] (EEE1) at (0,0) {};
\node[draw,circle,right=.8cm of EEE1] (EEE2) {};
\node[draw,circle,right=.8cm of EEE2] (EEE3) {};
\node[draw,circle,right=.8cm of EEE3] (EEE4) {};
\node[draw,circle,right=.8cm of EEE4] (EEE5) {};
\node[draw,circle,right=.8cm of EEE5] (EEE6) {};
\node[draw,circle,right=.8cm of EEE6] (EEE7) {};
\node[draw,circle,below=.8cm of EEE3] (EEEA) {};

\draw (EEE1) -- (EEE2)  node[] {};
\draw (EEE2) -- (EEE3)  node[] {};
\draw (EEE3) -- (EEE4)  node[] {};
\draw (EEE4) -- (EEE5)  node[] {};
\draw (EEE5) -- (EEE6)  node[] {};
\draw (EEE6) -- (EEE7)  node[] {};
\draw (EEE3) -- (EEEA)  node[] {};
\end{tikzpicture}
&&
{\footnotesize $\widetilde{E}_8$}

&
\begin{tikzpicture}[thick,scale=\scalevalue, every node/.style={transform shape}]
\node[draw,circle] (EEE1) at (0,0) {};
\node[draw,circle,right=.8cm of EEE1] (EEE2) {};
\node[draw,circle,right=.8cm of EEE2] (EEE3) {};
\node[draw,circle,right=.8cm of EEE3] (EEE4) {};
\node[draw,circle,right=.8cm of EEE4] (EEE5) {};
\node[draw,circle,right=.8cm of EEE5] (EEE6) {};
\node[draw,circle,right=.8cm of EEE6] (EEE7) {};
\node[draw,circle,right=.8cm of EEE7] (EEE8) {};
\node[draw,circle,below=.8cm of EEE3] (EEEA) {};

\draw (EEE1) -- (EEE2)  node[above,midway] {};
\draw (EEE2) -- (EEE3)  node[above,midway] {};
\draw (EEE3) -- (EEE4)  node[above,midway] {};
\draw (EEE4) -- (EEE5)  node[above,midway] {};
\draw (EEE5) -- (EEE6)  node[above,midway] {};
\draw (EEE6) -- (EEE7)  node[above,midway] {};
\draw (EEE8) -- (EEE7)  node[above,midway] {};
\draw (EEE3) -- (EEEA)  node[left,midway] {};
\end{tikzpicture}
\\
&&&&\\
\end{tabular}
\caption{The irreducible spherical Coxeter diagrams on the left and irreducible affine Coxeter diagrams on the right.}
\label{tab:classi_spherical_affine}
\end{table}

\section{Hyperbolic reflection groups}

\subsection{Hyperbolic polytopes}

Let $\R^{d,1}$ be the vector space $\R^{d+1}$ endowed with a non-degenerate symmetric bilinear form $\langle \cdot, \cdot \rangle$ of signature\footnote{The \emph{signature}\index{signature} of a symmetric matrix $B$ is the triple $(p,q,r)$ of the positive, negative, and zero indices of inertia of $B$. In the case $r=0$, we simply say that $B$ is of signature $(p,q)$.} $(d,1)$, and let $q$ be the associated quadratic form. A coordinate representation of $q$ with respect to some basis of $\R^{d,1}$ is:
$$ q(x) =  x_1^2 + \dotsm + x_d^2 - x_{d+1}^2.$$
A \emph{hyperbolic $d$-space}\index{hyperbolic space} $\Hb^{d}$ is a connected component of a hyperquadric: $$\Hb^{d} = \{ x \in  \R^{d,1} \mid q(x) = -1\ \textrm{ and }\, x_{d+1} > 0 \}.$$
The isometry group of $\Hb^{d}$ is $\mathrm{O}_{d,1}^+(\R)$, which consists of the elements of $\mathrm{O}_{d,1}(\R)$ that preserve $\Hb^{d}$. We often work with the projective model of $\Hb^{d}$: 
$$ \Hb^{d} = \{ x \in \R^{d,1} \mid q(x) < 0 \textrm{ and }\, x_{d+1} > 0 \}/{\R_+},$$
where the set $\R_+$ of positive scalars acts on $\R^{d,1} \smallsetminus \{ 0 \} $ by multiplication. If we set 
$$\mathbb{S}(\R^{d,1}) := {(\R^{d,1} \smallsetminus \{0\})}/{\R_+},$$ 
then $\Hb^{d}$ is an open subset of the \emph{projective sphere}\index{projective sphere} $\mathbb{S}(\R^{d,1})$. The closure $\overline{\Hb^{d}}$ of $\Hb^{d}$ in $\mathbb{S}(\R^{d,1})$ is the \emph{compactification}\index{hyperbolic space!compactification} of $\Hb^{d}$.

\medskip

A subset of $\Hb^d$ is a \emph{hyperbolic $d$-polytope}\index{hyperbolic polytope} if it is the intersection of a finite family of closed half-spaces of $\Hb^d$ and if it has non-empty interior. A \emph{hyperbolic Coxeter $d$-polytope}\index{hyperbolic Coxeter polytope} (or simply a $\Hb^d$-Coxeter polytope) is a hyperbolic $d$-polytope $P$ all its dihedral angles are sub-multiples of $\pi$. In other words, if two facets\footnote{A face of $P$ of codimension $1$ (resp. $2$) is called a \emph{facet}\index{facet} (resp. \emph{ridge}\index{ridge}) of $P$.}  $s$, $t$ of $P$ are adjacent,\footnote{Two facets $s$ and $t$ of $P$ are \emph{adjacent}\index{facet!adjacent} if $s \cap t$ is a ridge of $P$.} then the dihedral angle $\theta(s,t)$ between $s$ and $t$ is equal to $\nicefrac{\pi}{m}$ for some integer $m \geqslant 2$. When two facets $s, t$ are parallel, it is common to say that $\theta(s,t) = 0$.
A hyperbolic Coxeter polytope is \emph{right-angled}\index{hyperbolic Coxeter polytope!right-angled} if all its dihedral angles are $\nicefrac{\pi}{2}$.

\medskip

Associated with a hyperbolic Coxeter polytope $P$ is a Coxeter matrix $M = (m_{s,t})_{s,t \in S}$ on the set $S$ of facets of $P$: if $s, t \in S$ are adjacent, then $m_{s,t} = \nicefrac{\pi}{\theta(s,t)}$; otherwise, $m_{s,t} = \infty$. We denote by $W_P$ the Coxeter group of the Coxeter matrix $M$, and call it the \emph{Coxeter group of $P$}\index{hyperbolic Coxeter polytope!Coxeter group}.
If $s$ is a facet of $P$, then $\sigma_s$ denotes the reflection in the hyperplane containing $s$.

\begin{theorem}[{Poincaré \cite{Poincare1883}}]\label{thm:poincare}
Let $P$ be a $\Hb^d$-Coxeter polytope, and $W_P$ the Coxeter group of $P$. Then the homomorphism $\sigma : W_P \rightarrow \Isom (\Hb^d)$ defined by 
$$\sigma(s) = \sigma_{s} \quad\textrm{ for each } s \in S$$ is injective and the image $\G_P := \sigma(W_P)$ is discrete. Moreover, $P$ is a fundamental domain for the action of $\G_P$ on $\Hb^d$. In particular, if $P$ has finite volume (resp. is compact), then $\G_P$ is a lattice (resp. uniform lattice) of $\Isom (\Hb^d)$.
\end{theorem}

A subgroup $H$ of $\Isom (\Hb^d)$ is called a \emph{hyperbolic reflection group}\index{hyperbolic reflection group} if $H  = \G_P$ for some $\Hb^d$-Coxeter polytope $P$. In this case, we call $\G_P$ the \emph{reflection group of $P$}\index{hyperbolic Coxeter polytope!reflection group}. Theorem \ref{thm:poincare} provides a nice way to construct discrete subgroups of $\Isom (\Hb^d)$, even lattices. We shall review in the next section the classification of $\Hb^d$-Coxeter polytopes of finite volume in small dimensions and their non-existence in large dimensions.

\subsection{Classical results in dimensions 2 and 3}

\subsubsection{Dimension 2}

A necessary and sufficient condition for the existence of hyperbolic Coxeter $2$-polytopes of finite volume follows immediately from: 

\begin{theorem}
Let $\theta_1, \dotsc, \theta_n$ be real numbers such that $0 \leqslant \theta_i < \pi$ for each $i = 1, \dotsc, n$. Then there exists a hyperbolic polygon of finite volume with dihedral angles $\theta_1, \dotsc, \theta_n$ if and only if $ \sum_{i=1}^{n} \theta_i < (n-2)\pi$.
\end{theorem}

\subsubsection{Dimension 3}

Hyperbolic Coxeter $3$-polytopes of finite volume are well understood, notably thanks to the classification of hyperbolic $3$-polytopes with dihedral angles $\leqslant \nicefrac{\pi}{2}$ due to Andreev \cite{Andreev1,Andreev2}. During his PhD, Roeder found and fixed a gap in the original proof of Andreev (see \cite{Roeder_after_Andreev}). Hodgson and Rivin \cite{Hodgson_Rivin} gave a characterization of hyperbolic $3$-polytopes, which generalizes Andreev's theorem, in terms of a generalized Gauss map. 

\medskip
To express properly Andreev's theorem, one needs some definitions. Two compact polytopes $\mathcal{P}, \mathcal{P}'$ of the Euclidean space $\R^d$ are \emph{combinatorially equivalent}\index{combinatorially equivalent} if there is a bijection between their faces that preserves the inclusion relation. A combinatorial equivalence class is called a \emph{combinatorial polytope}\index{combinatorial polytope}. Note that if a hyperbolic polytope $P \subset \Hb^d$ is of finite volume, then the closure $\overline{P}$ of $P$ in $\overline{\Hb^d}$ is combinatorially equivalent to a compact polytope of $\R^d$. A \emph{labeled polytope}\index{labeled polytope} is a combinatorial polytope $\mathcal{P}$ with a labeling $\theta$, that is, a function from the ridges of $\mathcal{P}$ to $( 0 , \nicefrac{\pi}{2}]$. A hyperbolic polytope $P$ of finite volume \emph{realizes}\index{labeled polytope!realize} a labeled polytope $\mathcal{P}$ if there exists a combinatorial equivalence $\phi$ between the faces of $\overline{P}$ and $\mathcal{P}$ such that the dihedral angle at the ridge $e$ of $\overline{P}$ is the label $\theta( \phi(e) )$ of $\mathcal{P}$.

\medskip

Let $\mathcal{P}$ be a labeled $3$-polytope with labeling $\theta$. A \emph{$k$-circuit}\index{circuit} of $\mathcal{P}$ is a sequence of distinct facets $s_1, \ldots , s_k$ such that $e_i := s_i \cap s_{i+1}$ (indices are modulo $k$) is an edge of $\mathcal{P}$. A $k$-circuit is \emph{prismatic}\index{circuit!prismatic} if all the (closed) edges $e_i$ are disjoint. The \emph{angle sum}\index{circuit!angle sum} of a $k$-circuit is the real number $\sum_{i=1}^{k} \theta (e_i)$. A $k$-circuit is \emph{spherical}\index{circuit!spherical} (resp. \emph{Euclidean}\index{circuit!Euclidean}, resp. \emph{hyperbolic}\index{circuit!hyperbolic}) if its angle sum is bigger than (resp. equal to, resp. less than) $(k-2) \pi$. A vertex $v$ of $\mathcal{P}$ is \emph{spherical}\index{vertex!spherical} (resp. \emph{Euclidean}\index{vertex!Euclidean}) if the circuit consisting of the facets that contain $v$ is spherical (resp. Euclidean). The \emph{graph $W_\mathcal{P}$\index{labeled polytope!graph} of $\mathcal{P}$} is the graph whose nodes are the facets of $\mathcal{P}$ and such that two nodes $s,t$ are connected if and only if $s \cap t$ is not an edge, or $s \cap t$ is an edge and $\theta (s \cap t) < \nicefrac{\pi}{2}$.

\begin{theorem}[Andreev \cite{Andreev1,Andreev2}; see also Roeder--Hubbard--Dunbar \cite{Roeder_after_Andreev}]$\,$\\
Let $\mathcal{P}$ be a labeled $3$-polytope whose underlying polytope is not a tetrahedron. Then there exists a compact (resp. finite volume) hyperbolic $3$-polytope $P$ that realizes $\mathcal{P}$ if and only if:
\begin{enumerate}[label=(\roman*)]
	\item all the vertices of $\mathcal{P}$ are spherical (resp. spherical or Euclidean);
	\item all the prismatic $3$- and $4$-circuits of $\mathcal{P}$ are hyperbolic;
	\item the graph $W_\mathcal{P}$ of $\mathcal{P}$ is connected.
\end{enumerate}
In that case, the polytope $P$ is unique up to an isometry of $\Hb^3$.
\end{theorem}

\begin{remark}
A careful reader probably wonders what happen when the underlying polytope of $\mathcal{P}$ is a tetrahedron. This case is explained in \cite{Roeder_on_tetrahedra}. The list of all Coxeter tetrahedra of finite volume can be found in Tables \ref{Table_dim3_compact} and \ref{Table_dim3_noncompact}.
\end{remark}

\begin{remark}
The condition on the connectivity of $W_\mathcal{P}$ is often expressed by inequalities. One may notice that if the conditions $(i)$ and $(ii)$ are satisfied, then the disconnectedness of $W_\mathcal{P}$ implies that (see \cite[Prop. 1.5]{Roeder_after_Andreev}):
\begin{itemize}
	\item $\mathcal{P}$ is a right triangular prism, i.e., all the labels between the base facets and the joining facets are $\nicefrac{\pi}{2}$, or
	\item $\mathcal{P}$ is a quadrilateral pyramid with Euclidean apex such that the labels of two opposite edges in the base are $\nicefrac{\pi}{2}$.
\end{itemize}
\end{remark}

\subsection{The Gram matrix of a hyperbolic polytope}

A hyperbolic $d$-polytope $P$ is of the form
$$ P = \cap_{i=1}^{N} H_i^-,  $$
where $H_i^-$ is a closed half-space of $\Hb^d$ whose boundary is a hyperplane $H_i$. Each $H_i^-$ corresponds to a unique vector $u_i$ with the property that:
$$ H_i^- = \Hb^d \cap \{ x \in \R^{d,1} \mid  \langle x, u_i \rangle \leqslant 0 \}  \quad \textrm{and} \quad \langle u_i, u_i \rangle = 1.$$
The \emph{Gram matrix $G = (g_{i,j})$}\index{Gram matrix} of $P$ is a symmetric $N \times N$ matrix with entries $g_{i,j} = \langle u_i , u_j \rangle$.  A square matrix $B$ is \emph{reducible}\index{square matrix!reducible} if $B$ is the direct sum of smaller square matrices $B_1$ and $B_2$ (after a reordering of the indices), i.e., $B = \left(\begin{smallmatrix} B_1 & 0 \\ 0 & B_2\end{smallmatrix}\right)$. Otherwise, $B$ is \emph{irreducible}\index{square matrix!irreducible}. Every square matrix $B$ shall be the direct sum of irreducible submatrices, each of which we call a \emph{component}\index{square matrix!component} of $B$. The next theorem tells us that the hyperbolic polytopes are determined by the Gram matrices.

\begin{theorem}[Vinberg {\cite[Th.~2.1]{survey_hyp_reflec_vinberg}}]\label{thm:gram} Let $G = (g_{i,j})$ be an irreducible symmetric $N \times N$ matrix of signature $(d,1,N-d-1)$ such that the diagonal entries $g_{i,i} = 1$ and off-diagonal entries $g_{i,j}$ are $\leqslant 0$. Then there exists a hyperbolic $d$-polytope $P$ whose Gram matrix is $G$, and the polytope $P$ is unique up to an isometry of $\Hb^d$.
\end{theorem}

\begin{remark}
A detailed analysis of the Gram matrix can reveal the combinatorial structure of the hyperbolic polytope $P$ (see \cite[Th.~3.1 \& 3.2]{survey_hyp_reflec_vinberg}) and whether $P$ is compact or of finite volume (see \cite[Th.~4.1]{survey_hyp_reflec_vinberg}). 
\end{remark}

\subsection{Lannér and quasi-Lannér Coxeter groups}\label{subsec:Lanner}

A Coxeter group $W_S$ is \emph{Lannér}\index{Coxeter group!Lannér} (resp. \emph{quasi-Lannér}\index{Coxeter group!quasi-Lannér}) if $\det( \mathsf{C}_{W}) < 0$ and if for every proper subset $T \subset S$, the Coxeter group $W_T$ is spherical (resp. spherical or irreducible affine). If $P$ is a compact (resp. finite volume) Coxeter simplex, then the Coxeter group $W_P$ is Lannér (resp. quasi-Lannér). Conversely, if $W$ is a Lannér (resp. quasi-Lannér) Coxeter group of rank $d+1$, then its Cosine matrix $\mathsf{C}_W$ has signature $(d,1)$ and there exists a compact (resp. finite volume) $\Hb^d$-Coxeter polytope whose Gram matrix is $\frac{1}{2}\mathsf{C}_W$.

\medskip

In short, Lannér (resp. quasi-Lannér) Coxeter groups correspond to compact (resp. finite-volume) $\Hb^d$-Coxeter simplices. Such Coxeter simplices exist only in small dimensions.

\begin{theorem}[Lannér \cite{Lanner}, Koszul \cite{LectHypCoxGrKoszul} and Chein \cite{chein}]
If $W_S$ is a Lannér (resp. quasi-Lannér) Coxeter group, then $\# S \leqslant 5$ (resp. $\# S \leqslant 10$). Table \ref{tab:nb_quasi_lanner} indicates the number of Lannér (resp. quasi-Lannér) Coxeter groups of a given rank. The list of Lannér and quasi-Lannér Coxeter groups of rank $\geqslant 4$  can be found in \cite{chein}.

\newcommand{\siz}{0.15}
\newcommand{\size}{0.11}

\begin{table}[ht!]
\centering
\begin{tabular}{c|c|c}
\toprule[\siz em]
Dimension & $\sharp$ of quasi-Lannér & $\sharp$ of Lannér\\
$d = \#S - 1$     &  not Lannér                 &  Coxeter groups\\
\toprule[\siz em]
2 & $\infty$ & $\infty$ \\
\hline
3 & 23 & 9 \\
4 & 9 & 5\\
\hline
5 & 12& 0\\
6 &  3 & 0\\
7 &  4 & 0\\
8 &  4 & 0\\
9 & 3 & 0\\
\toprule[\siz em]
\end{tabular}
\vspace*{1em}
\caption{The numbers of quasi-Lannér or Lannér Coxeter groups}
\label{tab:nb_quasi_lanner}
\end{table}
\end{theorem}

\newcommand{\secondscalevalue}{0.50}

\begin{table}[ht!]
\centering
\begin{tabular}{cc cc cc}
\multicolumn{6}{c}{
\begin{tabular}{>{\centering\arraybackslash}m{.18\linewidth}>{\centering\arraybackslash}m{.18\linewidth}>{\centering\arraybackslash}m{.18\linewidth}>{\centering\arraybackslash}m{.18\linewidth}>{\centering\arraybackslash}m{.18\linewidth}}
\begin{tikzpicture}[thick,scale=\secondscalevalue, every node/.style={transform shape}]
\node[draw,circle] (C1) at (0,0) {};
\node[draw,circle,right=.8cm of C1] (C2) {};
\node[draw,circle,below=.8cm of C2] (C3) {}; 
\node[draw,circle,left=.8cm of C3] (C4) {};  
\draw (C1) -- (C2)  node[above,midway] {\large $4$};
\draw (C2) -- (C3)  node[above,midway] {};
\draw (C3) -- (C4) node[above,midway] {};
\draw (C4) -- (C1) node[above,midway] {};
\end{tikzpicture}
&
 \begin{tikzpicture}[thick,scale=\secondscalevalue, every node/.style={transform shape}]
\node[draw,circle] (D1) at (0,0) {};
\node[draw,circle,right=.8cm of D1] (D2) {};
\node[draw,circle,below=.8cm of D2] (D3) {}; 
\node[draw,circle,left=.8cm of D3] (D4) {};  
\draw (D1) -- (D2)  node[above,midway] {\large $5$};
\draw (D2) -- (D3)  node[above,midway] {};
\draw (D3) -- (D4) node[above,midway] {};
\draw (D4) -- (D1) node[above,midway] {};
\end{tikzpicture}
&
 \begin{tikzpicture}[thick,scale=\secondscalevalue, every node/.style={transform shape}]
\node[draw,circle] (E1) at (0,0) {};
\node[draw,circle,right=.8cm of E1] (E2) {};
\node[draw,circle,below=.8cm of E2] (E3) {}; 
\node[draw,circle,left=.8cm of E3] (E4) {};
\draw (E1) -- (E2)  node[above,midway] {\large $4$};
\draw (E2) -- (E3)  node[above,midway] {};
\draw (E3) -- (E4) node[below,midway] {\large $4$};
\draw (E4) -- (E1) node[above,midway] {};
\end{tikzpicture}
 &
\begin{tikzpicture}[thick,scale=\secondscalevalue, every node/.style={transform shape}]
\node[draw,circle] (F1) at (0,0) {};
\node[draw,circle,right=.8cm of F1] (F2) {};
\node[draw,circle,below=.8cm of F2] (F3) {}; 
\node[draw,circle,left=.8cm of F3] (F4) {};
\draw (F1) -- (F2)  node[above,midway] {\large $5$};
\draw (F2) -- (F3)  node[above,midway] {};
\draw (F3) -- (F4) node[below,midway] {\large $4$};
\draw (F4) -- (F1) node[above,midway] {};
\end{tikzpicture}
 &
 \begin{tikzpicture}[thick,scale=\secondscalevalue, every node/.style={transform shape}]
\node[draw,circle] (G1) at (0,0) {};
\node[draw,circle,right=.8cm of G1] (G2) {};
\node[draw,circle,below=.8cm of G2] (G3) {}; 
\node[draw,circle,left=.8cm of G3] (G4) {};
  \draw (G1) -- (G2)  node[above,midway] {\large $5$};
\draw (G2) -- (G3)  node[above,midway] {};
\draw (G3) -- (G4) node[below,midway] {\large $5$};
\draw (G4) -- (G1) node[above,midway] {};
\end{tikzpicture}
\\
\end{tabular}
}
\\
\\
\multicolumn{6}{c}{
\begin{tabular}{>{\centering\arraybackslash}m{.22\linewidth}>{\centering\arraybackslash}m{.22\linewidth}>{\centering\arraybackslash}m{.22\linewidth}>{\centering\arraybackslash}m{.22\linewidth}}
\begin{tikzpicture}[thick,scale=\secondscalevalue, every node/.style={transform shape}]
\node[draw,circle] (C1) at (0,0) {};
 \node[draw,circle,right=.8cm of C1] (C2) {};
\node[draw,circle,right=.8cm of C2] (C3) {}; 
\node[draw,circle,right=.8cm of C3] (C4) {};
\draw (C1) -- (C2)  node[above,midway] {};
\draw (C2) -- (C3)  node[above,midway] {\large $5$};
\draw (C3) -- (C4) node[] {};
\end{tikzpicture}
&
\begin{tikzpicture}[thick,scale=\secondscalevalue, every node/.style={transform shape}]
\node[draw,circle] (CC1) at (0,0) {};
\node[draw,circle,right=.8cm of CC1] (CC2) {};
\node[draw,circle,right=.8cm of CC2] (CC3) {}; 
\node[draw,circle,right=1cm of CC3] (CC4) {};
\draw (CC1) -- (CC2)  node[above,midway] {\large $5$};
\draw (CC2) -- (CC3)  node[above,midway] {};
\draw (CC3) -- (CC4) node[above,midway] {\large $4$};
\end{tikzpicture}
 &
\begin{tikzpicture}[thick,scale=\secondscalevalue, every node/.style={transform shape}]
\node[draw,circle] (D1) at (0,0) {};
\node[draw,circle,right=.8cm of D1] (D2) {};
\node[draw,circle,right=.8cm of D2] (D3) {}; 
\node[draw,circle,right=1cm of D3] (D4) {};
\draw (D1) -- (D2)  node[above,midway] {\large $5$};
\draw (D2) -- (D3)  node[above,midway] {};
\draw (D3) -- (D4) node[above,midway] {\large $5$};
\end{tikzpicture}
 &
\begin{tikzpicture}[thick,scale=\secondscalevalue, every node/.style={transform shape}]
\node[draw,circle] (B4) at (0,0) {};
\node[draw,circle,right=.8cm of B4] (B5) {};
\node[draw,circle,above right=.8cm of B5] (B6) {};
\node[draw,circle,below right=.8cm of B5] (B7) {};
\draw (B4) -- (B5) node[above,midway] {\large $5$};
\draw (B5) -- (B6) node[above,midway] {};
\draw (B5) -- (B7) node[above,midway] {};
\end{tikzpicture}
\end{tabular}
}
\end{tabular}
\caption{The nine compact hyperbolic tetrahedra}
\label{Table_dim3_compact}
\end{table}

\begin{table}[ht!]
\centering
\begin{tabular}{cc cc cc}
\multicolumn{6}{c}{
\begin{tabular}{>{\centering\arraybackslash}m{.18\linewidth}>{\centering\arraybackslash}m{.18\linewidth}>{\centering\arraybackslash}m{.18\linewidth}>{\centering\arraybackslash}m{.18\linewidth}>{\centering\arraybackslash}m{.18\linewidth}}
\begin{tikzpicture}[thick,scale=\secondscalevalue, every node/.style={transform shape}]
\node[draw,circle] (C1) at (0,0) {};
\node[draw,circle,right=.8cm of C1] (C2) {};
\node[draw,circle,below=.8cm of C2] (C3) {}; 
\node[draw,circle,left=.8cm of C3] (C4) {};
\draw (C1) -- (C2)  node[above,midway] {\large $4$};
\draw (C2) -- (C3)  node[above,midway] {};
\draw (C3) -- (C4) node[above,midway] {};
\draw (C4) -- (C1) node[left,midway] {\large $4$};
\end{tikzpicture}
&
 \begin{tikzpicture}[thick,scale=\secondscalevalue, every node/.style={transform shape}]
\node[draw,circle] (D1) at (0,0) {};
\node[draw,circle,right=.8cm of D1] (D2) {};
\node[draw,circle,below=.8cm of D2] (D3) {}; 
\node[draw,circle,left=.8cm of D3] (D4) {};
\draw (D1) -- (D2)  node[above,midway] {\large $4$};
\draw (D2) -- (D3)  node[right,midway] {\large $4$};
\draw (D3) -- (D4) node[above,midway] {};
\draw (D4) -- (D1) node[left,midway] {\large $4$};
\end{tikzpicture}
&
 \begin{tikzpicture}[thick,scale=\secondscalevalue, every node/.style={transform shape}]
\node[draw,circle] (E1) at (0,0) {};
\node[draw,circle,right=.8cm of E1] (E2) {};
\node[draw,circle,below=.8cm of E2] (E3) {}; 
\node[draw,circle,left=.8cm of E3] (E4) {};
\draw (D1) -- (D2)  node[above,midway] {\large $4$};
\draw (D2) -- (D3)  node[right,midway] {\large $4$};
\draw (D3) -- (D4) node[below,midway] {\large $4$};
\draw (D4) -- (D1) node[left,midway] {\large $4$};
\end{tikzpicture}
 &
 \begin{tikzpicture}[thick,scale=\secondscalevalue, every node/.style={transform shape}]
\node[draw,circle] (F1) at (0,0) {};
\node[draw,circle,right=.8cm of F1] (F2) {};
\node[draw,circle,below=.8cm of F2] (F3) {}; 
\node[draw,circle,left=.8cm of F3] (F4) {};
\draw (F1) -- (F2)  node[above,midway] {\large $6$};
\draw (F2) -- (F3)  node[above,midway] {};
\draw (F3) -- (F4) node[below,midway] {};
\draw (F4) -- (F1) node[above,midway] {};
\end{tikzpicture}
 &
 \begin{tikzpicture}[thick,scale=\secondscalevalue, every node/.style={transform shape}]
\node[draw,circle] (F1) at (0,0) {};
\node[draw,circle,right=.8cm of F1] (F2) {};
\node[draw,circle,below=.8cm of F2] (F3) {}; 
\node[draw,circle,left=.8cm of F3] (F4) {};
\draw (F1) -- (F2)  node[above,midway] {\large $6$};
\draw (F2) -- (F3)  node[above,midway] {};
\draw (F3) -- (F4) node[below,midway] {\large $4$};
\draw (F4) -- (F1) node[above,midway] {};
\end{tikzpicture}
\\
\end{tabular}
}
\\
\\
\multicolumn{6}{c}{
\begin{tabular}{>{\centering\arraybackslash}m{.18\linewidth}>{\centering\arraybackslash}m{.18\linewidth}>{\centering\arraybackslash}m{.18\linewidth}>{\centering\arraybackslash}m{.18\linewidth}>{\centering\arraybackslash}m{.18\linewidth}}
 \begin{tikzpicture}[thick,scale=\secondscalevalue, every node/.style={transform shape}]
\node[draw,circle] (F1) at (0,0) {};
\node[draw,circle,right=.8cm of F1] (F2) {};
\node[draw,circle,below=.8cm of F2] (F3) {}; 
\node[draw,circle,left=.8cm of F3] (F4) {};
\draw (F1) -- (F2)  node[above,midway] {\large $6$};
\draw (F2) -- (F3)  node[above,midway] {};
\draw (F3) -- (F4) node[below,midway] {\large $5$};
\draw (F4) -- (F1) node[above,midway] {};
\end{tikzpicture}
&
 \begin{tikzpicture}[thick,scale=\secondscalevalue, every node/.style={transform shape}]
\node[draw,circle] (G1) at (0,0) {};
\node[draw,circle,right=.8cm of G1] (G2) {};
\node[draw,circle,below=.8cm of G2] (G3) {}; 
\node[draw,circle,left=.8cm of G3] (G4) {};
\draw (G1) -- (G2)  node[above,midway] {\large $6$};
\draw (G2) -- (G3)  node[above,midway] {};
\draw (G3) -- (G4) node[below,midway] {\large $6$};
\draw (G4) -- (G1) node[above,midway] {};
\end{tikzpicture}
&
 \begin{tikzpicture}[thick,scale=\secondscalevalue, every node/.style={transform shape}]
\node[draw,circle] (G1) at (0,0) {};
\node[draw,circle,right=.8cm of G1] (G2) {};
\node[draw,circle,below=.8cm of G2] (G3) {}; 
\node[draw,circle,left=.8cm of G3] (G4) {};
\draw (G1) -- (G2)  node[above,midway] {};
\draw (G2) -- (G3)  node[above,midway] {};
\draw (G3) -- (G4) node[below,midway] {};
\draw (G4) -- (G1) node[above,midway] {};
\draw (G2) -- (G4) node[above,midway] {};
\end{tikzpicture}
 &
 &
\\
\end{tabular}
}
\\
\\
\multicolumn{6}{c}{
\begin{tabular}{>{\centering\arraybackslash}m{.18\linewidth}>{\centering\arraybackslash}m{.18\linewidth}>{\centering\arraybackslash}m{.18\linewidth}>{\centering\arraybackslash}m{.18\linewidth}>{\centering\arraybackslash}m{.18\linewidth}}
\begin{tikzpicture}[thick,scale=\secondscalevalue, every node/.style={transform shape}]
\node[draw,circle] (C1) at (0,0) {};
 \node[draw,circle] (C2) at (-180:1) {};
 \node[draw,circle] (C3) at (60:1) {};
\node[draw,circle] (C4) at (-60:1) {};
 \draw (C1) -- (C2)  node[above,midway] {};
\draw (C3) -- (C1)  node[above,midway] {};
\draw (C4) -- (C1) node[left,midway] {};
 \draw (C3) -- (C2)  node[above,midway] {};
\draw (C3) -- (C4)  node[above,midway] {};
\draw (C4) -- (C2) node[left,midway] {};
\end{tikzpicture}
&
\begin{tikzpicture}[thick,scale=\secondscalevalue, every node/.style={transform shape}]
\node[draw,circle] (C1) at (-.8,0) {};
 \node[draw,circle] (C2) at (60:.8) {};
 \node[draw,circle] (C3) at (-60:.8) {};
\node[draw,circle] (C4) at (-2,0) {};
  \draw (C1) -- (C2)  node[above,midway] {};
\draw (C2) -- (C3)  node[above,midway] {};
\draw (C3) -- (C1) node[above,midway] {};
\draw (C4) -- (C1) node[above,midway] {};
\end{tikzpicture}
&
\begin{tikzpicture}[thick,scale=\secondscalevalue, every node/.style={transform shape}]
\node[draw,circle] (C1) at (-.8,0) {};
 \node[draw,circle] (C2) at (60:.8) {};
 \node[draw,circle] (C3) at (-60:.8) {};
\node[draw,circle] (C4) at (-2,0) {};
  \draw (C1) -- (C2)  node[above,midway] {};
\draw (C2) -- (C3)  node[above,midway] {};
\draw (C3) -- (C1) node[above,midway] {};
\draw (C4) -- (C1) node[above,midway] {\large $4$};
\end{tikzpicture}
 &
\begin{tikzpicture}[thick,scale=\secondscalevalue, every node/.style={transform shape}]
\node[draw,circle] (C1) at (-.8,0) {};
 \node[draw,circle] (C2) at (60:.8) {};
 \node[draw,circle] (C3) at (-60:.8) {};
\node[draw,circle] (C4) at (-2,0) {};
  \draw (C1) -- (C2)  node[above,midway] {};
\draw (C2) -- (C3)  node[above,midway] {};
\draw (C3) -- (C1) node[above,midway] {};
\draw (C4) -- (C1) node[above,midway] {\large $5$};
\end{tikzpicture}
 &
\begin{tikzpicture}[thick,scale=\secondscalevalue, every node/.style={transform shape}]
\node[draw,circle] (C1) at (-.8,0) {};
 \node[draw,circle] (C2) at (60:.8) {};
 \node[draw,circle] (C3) at (-60:.8) {};
\node[draw,circle] (C4) at (-2,0) {};
  \draw (C1) -- (C2)  node[above,midway] {};
\draw (C2) -- (C3)  node[above,midway] {};
\draw (C3) -- (C1) node[above,midway] {};
\draw (C4) -- (C1) node[above,midway] {\large $6$};
\end{tikzpicture}
\\
\end{tabular}
}
\\
\\
\multicolumn{6}{c}{
\begin{tabular}{>{\centering\arraybackslash}m{.22\linewidth}>{\centering\arraybackslash}m{.22\linewidth}>{\centering\arraybackslash}m{.22\linewidth}>{\centering\arraybackslash}m{.22\linewidth}}
\begin{tikzpicture}[thick,scale=\secondscalevalue, every node/.style={transform shape}]
\node[draw,circle] (C1) at (0,0) {};
 \node[draw,circle,right=.8cm of C1] (C2) {};
\node[draw,circle,right=.8cm of C2] (C3) {}; 
\node[draw,circle,right=.8cm of C3] (C4) {};
  \draw (C1) -- (C2)  node[above,midway] {};
\draw (C2) -- (C3)  node[above,midway] {\large $6$};
\draw (C3) -- (C4) node[above,midway] {};
\end{tikzpicture}
&
\begin{tikzpicture}[thick,scale=\secondscalevalue, every node/.style={transform shape}]
\node[draw,circle] (C1) at (0,0) {};
 \node[draw,circle,right=.8cm of C1] (C2) {};
\node[draw,circle,right=.8cm of C2] (C3) {}; 
\node[draw,circle,right=.8cm of C3] (C4) {};
\draw (C1) -- (C2)  node[above,midway] {\large $4$};
\draw (C2) -- (C3)  node[above,midway] {\large $4$};
\draw (C3) -- (C4) node[above,midway] {};
\end{tikzpicture}
 &
\begin{tikzpicture}[thick,scale=\secondscalevalue, every node/.style={transform shape}]
\node[draw,circle] (C1) at (0,0) {};
 \node[draw,circle,right=.8cm of C1] (C2) {};
\node[draw,circle,right=.8cm of C2] (C3) {}; 
\node[draw,circle,right=.8cm of C3] (C4) {};
  \draw (C1) -- (C2)  node[above,midway] {\large $4$};
\draw (C2) -- (C3)  node[above,midway] {\large $4$};
\draw (C3) -- (C4) node[above,midway] {\large $4$};
\end{tikzpicture}
 &
\begin{tikzpicture}[thick,scale=\secondscalevalue, every node/.style={transform shape}]
\node[draw,circle] (C1) at (0,0) {};
 \node[draw,circle,right=.8cm of C1] (C2) {};
\node[draw,circle,right=.8cm of C2] (C3) {}; 
\node[draw,circle,right=.8cm of C3] (C4) {};
 \draw (C1) -- (C2)  node[above,midway] {\large $6$};
\draw (C2) -- (C3)  node[above,midway] {};
\draw (C3) -- (C4) node[above,midway] {};
\end{tikzpicture}
\end{tabular}
}
\\
\\
\multicolumn{6}{c}{
\begin{tabular}{>{\centering\arraybackslash}m{.22\linewidth}>{\centering\arraybackslash}m{.22\linewidth}>{\centering\arraybackslash}m{.22\linewidth}>{\centering\arraybackslash}m{.22\linewidth}}
\begin{tikzpicture}[thick,scale=\secondscalevalue, every node/.style={transform shape}]
\node[draw,circle] (C1) at (0,0) {};
 \node[draw,circle,right=.8cm of C1] (C2) {};
\node[draw,circle,right=.8cm of C2] (C3) {}; 
\node[draw,circle,right=.8cm of C3] (C4) {};
 \draw (C1) -- (C2)  node[above,midway] {\large $6$};
\draw (C2) -- (C3)  node[above,midway] {};
\draw (C3) -- (C4) node[above,midway] {\large $4$};
\end{tikzpicture}
&
\begin{tikzpicture}[thick,scale=\secondscalevalue, every node/.style={transform shape}]
\node[draw,circle] (C1) at (0,0) {};
 \node[draw,circle,right=.8cm of C1] (C2) {};
\node[draw,circle,right=.8cm of C2] (C3) {}; 
\node[draw,circle,right=.8cm of C3] (C4) {};
 \draw (C1) -- (C2)  node[above,midway] {\large $6$};
\draw (C2) -- (C3)  node[above,midway] {};
\draw (C3) -- (C4) node[above,midway] {\large $5$};
\end{tikzpicture}
 &
\begin{tikzpicture}[thick,scale=\secondscalevalue, every node/.style={transform shape}]
\node[draw,circle] (C1) at (0,0) {};
 \node[draw,circle,right=.8cm of C1] (C2) {};
\node[draw,circle,right=.8cm of C2] (C3) {}; 
\node[draw,circle,right=.8cm of C3] (C4) {};
  \draw (C1) -- (C2)  node[above,midway] {\large $6$};
\draw (C2) -- (C3)  node[above,midway] {};
\draw (C3) -- (C4) node[above,midway] {\large $6$};
\end{tikzpicture}
 &
\end{tabular}
}
\\
\\
\multicolumn{6}{c}{
\begin{tabular}{>{\centering\arraybackslash}m{.22\linewidth}>{\centering\arraybackslash}m{.22\linewidth}>{\centering\arraybackslash}m{.22\linewidth}>{\centering\arraybackslash}m{.22\linewidth}}
\begin{tikzpicture}[thick,scale=\secondscalevalue, every node/.style={transform shape}]
\node[draw,circle] (C1) at (0,0) {};
 \node[draw,circle,right=.8cm of C1] (C2) {};
\node[draw,circle,right=.8cm of C2] (C3) {}; 
\node[draw,circle,below=.8cm of C2] (C4) {};
 \draw (C1) -- (C2)  node[above,midway] {};
\draw (C2) -- (C3)  node[above,midway] {\large $6$};
\draw (C2) -- (C4) node[above,midway] {};
\end{tikzpicture}
&
\begin{tikzpicture}[thick,scale=\secondscalevalue, every node/.style={transform shape}]
\node[draw,circle] (C1) at (0,0) {};
 \node[draw,circle,right=.8cm of C1] (C2) {};
\node[draw,circle,right=.8cm of C2] (C3) {}; 
\node[draw,circle,below=.8cm of C2] (C4) {};
 \draw (C1) -- (C2)  node[above,midway] {\large $4$};
\draw (C2) -- (C3)  node[above,midway] {\large $4$};
\draw (C2) -- (C4) node[above,midway] {};
\end{tikzpicture}
 &
\begin{tikzpicture}[thick,scale=\secondscalevalue, every node/.style={transform shape}]
\node[draw,circle] (C1) at (0,0) {};
 \node[draw,circle,right=.8cm of C1] (C2) {};
\node[draw,circle,right=.8cm of C2] (C3) {}; 
\node[draw,circle,below=.8cm of C2] (C4) {};
 \draw (C1) -- (C2)  node[above,midway] {\large $4$};
\draw (C2) -- (C3)  node[above,midway] {\large $4$};
\draw (C2) -- (C4) node[left,midway] {\large $4$};
\end{tikzpicture}
 &
\end{tabular}
}
\end{tabular}
\caption{The twenty-three hyperbolic tetrahedra of finite volume which are not compact}
\label{Table_dim3_noncompact}
\end{table}

The non-existence of hyperbolic Coxeter simplices in large dimensions is in fact the first side of a more general non-existence theorem in Section \ref{subsec:absence}.

\begin{remark}
Tables \ref{Table_dim3_compact} and \ref{Table_dim3_noncompact} give us the list of all the Lannér or quasi-Lannér Coxeter groups of rank $4$, which correspond to compact or finite volume hyperbolic tetrahedra, respectively.
\end{remark}

\subsection{Absence in large dimension}\label{subsec:absence}

\begin{theorem}[Vinberg \cite{absence_vinberg} and Prokhorov \cite{absence_prokhorov}]\label{thm:absence}
If $\Gamma_P$ is a discrete reflection group of $\mathrm{Isom}(\HH^d)$ with compact (resp. finite volume) fundamental domain $P$, then $d \leqslant 29$ (resp. $d \leqslant 995$).
\end{theorem}

The proof of Vinberg (resp. Prokhorov) uses Nikulin's inequality for simple polytopes (resp. edge-simple polytopes) established in \cite{lemma_nikulin} (resp. \cite{lemma_khovanskii}). The upper bounds in the right-angled case are better:

\begin{theorem}[Potyagailo--Vinberg \cite{potyagailo_vinberg} and Dufour \cite{dufour}]
If $\Gamma_P$ is a discrete reflection group of $\mathrm{Isom}  (\HH^d)$ with compact (resp. finite volume) right-angled fundamental domain $P$, then $d \leqslant 4$  (resp. $d \leqslant 12$).
\end{theorem}


Except for compact right-angled polytopes, the upper bounds are far from being sharp.

\begin{itemize}
	\item There exists a compact right-angled hyperbolic 4-polytope: 120-cell.
	
	\item Examples of finite volume right-angled $d$-polytopes are known in dimension $d \leqslant 8$ (see \cite{potyagailo_vinberg}). 
	
	\item Examples of compact $d$-polytopes are known in dimension $d \leqslant  8$ (see \cite{bugaenko1,bugaenko_world_record} for $d=7,8$).

	\item Examples of finite volume $d$-polytopes are known in dimension $ d \leqslant 21$ and $d \neq 20$ (see \cite{allcock} and the references therein for $d \leqslant 19$ and \cite{borcherds} for $d=21$).
\end{itemize}

\begin{remark}
Moussong observed that the argument of Vinberg \cite{absence_vinberg} may be extended to show that if a Coxeter group $W$ is word hyperbolic and the nerve of $W$ is a generalized homology $(d-1)$-sphere, then $d \leqslant 29$ (see \cite[Prop.~12.6.7]{davis_book}). Here, the \emph{nerve}\index{Coxeter group!nerve} of a Coxeter group $W_S$ is the poset of all nonempty spherical subsets of $S$ partially ordered by inclusion, which is an abstract simplicial complex, and a \emph{generalized homology $d$-sphere}\index{generalized homology sphere} is a homology $d$-manifold with the same homology as the $d$-sphere.
\end{remark}

\subsection{Hyperbolic Coxeter polytopes with few facets}

The complete classification of compact or finite-volume hyperbolic polytopes is not an easy task. Only compact $d$-polytopes with $N$ facets $(N \leqslant d+3)$ and finite-volume $d$-polytopes with $N$ facets $(N \leqslant d+2)$ were classified. For more details, we refer the reader to the web page maintained by Felikson and Tumarkin:

\href{https://www.maths.dur.ac.uk/users/anna.felikson/Polytopes/polytopes.html}{https://www.maths.dur.ac.uk/users/anna.felikson/Polytopes/polytopes.html}

\noindent
This webpage contains all the known examples of hyperbolic Coxeter polytopes of dimension $\geqslant 4$.

\subsection{Convex cocompact hyperbolic reflection groups}

A subgroup $\Gamma$ of $\mathrm{Isom}(\mathbb{H}^d)$ is \emph{convex cocompact}\index{convex cocompact} if there exists a $\Gamma$-invariant convex subset $\mathscr{C}$ of $\mathbb{H}^d$ such that $\Gamma$ acts properly discontinuously on $\mathscr{C}$ with compact quotient. Using the following theorems, one can easily check when a reflection group $\Gamma_P$ of $\mathrm{Isom}(\mathbb{H}^d)$ is convex cocompact.

\begin{theorem}[Desgroseilliers--Haglund {\cite[Th.~4.12]{Desgroseilliers_Haglund}}]\label{thm:DH}
Let $P$ be an $\Hb^d$-Coxeter polytope and let $\Gamma_P$ be its reflection group. Then $\Gamma_P$ is convex cocompact if and only if 
\begin{enumerate}[label=(\roman*)]
\item $\Gamma_P$ is word-hyperbolic, and
\item $P$ has no pair of asymptotic facets.
\end{enumerate}  
\end{theorem}

\begin{remark}
The condition $(ii)$ may be replaced by $(ii')$ there is no pair of facets $s, t$ of $P$ such that $g_{s,t} = -1$, where $G = (g_{s,t})$ is the Gram matrix of $P$. 
\end{remark}

\begin{theorem}[{Moussong's hyperbolicity criterion \cite{moussong}}]\label{thm:moussong} Let $W_S$ be a Coxeter group. Then $W_S$ is word hyperbolic if and only if $S$ does not contain two orthogonal non-spherical subsets, nor any affine subset of rank $\geqslant 3$.
\end{theorem}

\begin{remark}
Desgroseilliers and Haglund \cite[Th.~1.1]{Desgroseilliers_Haglund} found a class of Coxeter groups which can be realized as a convex cocompact subgroup of $\mathrm{Isom}(\mathbb{H}^d)$, which is not a reflection group, and they conjectured that there exists a word-hyperbolic Coxeter group which admits a convex cocompact representation into $\mathrm{Isom}(\mathbb{H}^d)$ but which cannot be realized as a convex cocompact \emph{reflection} group of $\mathrm{Isom}(\mathbb{H}^n)$ for any $n \in \mathbb{N}$.
\end{remark}

\section{Projective reflection groups}

\subsection{Tits--Vinberg's Theorem}

Let $V$ be a vector space over $\R$, and let $\PS(V)$ be the projective sphere. We denote by $\SL^{\pm}(V)$ the group of automorphism of $\PS(V)$, i.e., $$\SL^{\pm}(V) = \{ g \in \GL(V) \mid \det(g) = \pm 1 \}.$$ 
We denote by $\hat{\PS}$ the natural projection of $V \smallsetminus \{ 0\}$ to $\PS(V)$, and let $\PS(W) := \hat{\PS}(W  \smallsetminus \{ 0 \} )$ for any subset $W$ of $V$.
The complement of a projective hyperplane in $\PS(V)$ consists of two connected components, each of which we call an \emph{affine chart}\index{affine chart} of $\PS(V)$. A \emph{cone}\index{cone} is a subset of $V$ which is invariant under multiplication by positive scalars.
A subset $\mathscr{C}$ of $\PS(V)$ is \emph{convex}\index{convex} if there exists a convex cone $U$ of $V$ such that $\mathscr{C} = \PS(V)$, \emph{properly convex}\index{properly convex} if it is convex and its closure lies in some affine chart, and \emph{strictly convex}\index{strictly convex} if in addition its boundary does not contain any nontrivial projective line segment. Hyperbolic spaces are special examples of strictly convex open subsets of $\PS(V)$.

\medskip

A \emph{projective polytope}\index{projective polytope} is a properly convex subset $P$ of $\PS(V)$ such that $P$ has a non-empty interior and $P = \cap_{i=1}^{N} \PS( \{x \in V \mid \alpha_i(x) \leqslant 0 \} )$, where $\alpha_i$, $i=1, \dotsc, N$, are linear forms on $V$. We always assume that $P$ has $N$ facets, i.e., to define $P$, we need all the $N$ linear forms $(\alpha_i)_{i=1}^{N}$. A \emph{projective reflection}\index{projective reflection} is an element of $\SL^{\pm}(V)$ of order $2$ which is the identity on a hyperplane. Every projective reflection $\sigma$ can be written as:
$$ \sigma = \mathrm{Id} - \alpha \otimes v, \quad \textrm{i.e.,} \quad \sigma( x ) =  x - \alpha( x ) v\quad \forall x \in V, $$
where $\alpha$ is a linear form on $V$ and $v$ is a vector of $V$ such that $\alpha(v)=2$.

\medskip

Let $P$ be a projective polytope and $S$ the set of facets of $P$. A \emph{reflection in a facet} $s \in S$ is a projective reflection $\sigma_s$ which fixes each point of $s$. A \emph{pre-mirror polytope}\index{pre-mirror polytope} is a projective polytope $P$ together with one reflection $\sigma_s$ in each facet $s$ of $P$. So, one may choose $\sigma_s = \mathrm{Id} - \alpha_s \otimes v_s$ with $\alpha_s(v_s)=2$ such that
$P = \cap_{s \in S} \PS( \{ x \in V \mid \alpha_s(x) \leqslant 0 \} )$. Note that the pairs $(\alpha_s,v_s)$ are uniquely determined only up to multiplication by a positive real number.

\medskip

If $P$ is a pre-mirror polytope, then $\G_P$ denotes the group generated by the reflections in the facets of $P$. We say that $\G_P$ is a \emph{projective reflection group}\index{projective reflection group} if for any $\g \in \G_P$,
$$
\g(\mathring{P}) \cap \mathring{P} \neq \varnothing
 \quad \Rightarrow \quad
\gamma = \mathrm{Id},
$$
where $\mathring{P}$ denotes the interior of $P$. 

In the next paragraph, we introduce a relevant tool to formulate Proposition \ref{prop:necessary} and Theorem \ref{thm:Tits_Vinberg} which express necessary and sufficient conditions for $\Gamma_P$ of a pre-mirror polytope $P$ to be a projective reflection group. A key notion is that of Cartan matrix of mirror polytope which generalizes the twice of the Gram matrix of hyperbolic polytope.
\begin{definition}\label{def:cartan_matrix}
A \emph{Cartan matrix}\index{Cartan matrix} on a set $S$ is a $S \times S$ matrix $\Cart_S = (a_{s,t})_{s,t \in S}$ which satisfies the conditions: $(i)$ $a_{s,s} = 2,\, \forall s \in S$; $(ii)$ $a_{s,t} \leqslant 0,\, \forall s \neq t \in S$; $(iii)$ $a_{s,t} = 0 \Leftrightarrow a_{t,s} = 0,\, \forall s\neq t \in S$.
\end{definition}

\medskip

A \emph{mirror polytope}\index{mirror polytope} is a pre-mirror polytope $P$ such that the matrix $\Cart_P := (\alpha_s(v_t))_{s,t \in S}$ is a Cartan matrix. In this case, we call $\Cart_P$ the \emph{Cartan matrix of $P$}\index{mirror polytope!Cartan matrix}. A Cartan matrix $\Cart$ on $S$ is of \emph{Coxeter type}\index{Cartan matrix!Coxeter type} if for any $s \neq t \in S$,
$$
a_{s,t} a_{t,s} < 4 \quad
\Rightarrow \quad
\frac{\pi}{\arccos \left( \frac{1}{2} \sqrt{a_{s,t} a_{t,s}}  \right)} \in \N.
$$
A \emph{projective Coxeter polytope}\index{projective Coxeter polytope} is a mirror polytope $P$ whose Cartan matrix $\Cart_P$ is of Coxeter type. For each pair of adjacent facets $s, t$ of $P$, the \emph{dihedral angle}\index{dihedral angle} of the ridge $s \cap t$ is said to be $\nicefrac{\pi}{m_{s,t}}$ if $a_{s,t} a_{t,s} = 4 \cos^2 (\nicefrac{\pi}{m_{s,t}})$.
 
\begin{proposition}[Vinberg {\cite[Prop.~17]{bible}}]\label{prop:necessary}
Let $P$ be a pre-mirror polytope. If the group $\G_P$ is a projective reflection group, then $P$ is a projective Coxeter polytope.
\end{proposition}

A Cartan matrix $\Cart_S$ and a Coxeter group $W_S$ are \emph{compatible}\index{Cartan matrix!compatible} when:
\begin{enumerate}[label=(\roman*)]
	\item $\forall s,t \in S,\quad  m_{s,t} = 2 \,\, \Leftrightarrow \,\, a_{s,t} = 0$;
	\item $\forall s,t \in S,\quad m_{s,t} < \infty \,\, \Leftrightarrow \,\, a_{s,t} a_{t,s} = 4 \cos^2 ( \nicefrac{\pi}{m_{s,t}} )$;
	\item $\forall s,t \in S,\quad m_{s,t} = \infty \,\, \Leftrightarrow \,\, a_{s,t} a_{t,s} \geqslant 4$.
\end{enumerate}
\bigskip

It is clear that there is at most one Coxeter group compatible with a given Cartan matrix, and that a Cartan matrix $\Cart_S$ is compatible with some Coxeter group $W_S$ if and only if $\Cart_S$ is of Coxeter type. If $P$ is a projective Coxeter polytope, then $W_P$ denotes the unique Coxeter group compatible with $\Cart_P$. The following is a generalization of Theorem \ref{thm:poincare} to the projective setting.

\begin{theorem}[Bourbaki {\cite[Chap.~V]{Bourbaki_group_456} and Vinberg \cite[Th.~2]{bible}
 }]\label{thm:Tits_Vinberg}\footnote{Theorem \ref{thm:Tits_Vinberg} was proved by Tits for $\Delta_W$, which we define in Section \ref{subsec:CartanToPolytope}, and by Vinberg for the general case.}
Let $P$ be a projective Coxeter polytope of $\PS(V)$ with Coxeter group $W_P$, and let $\Gamma_P$ be the group generated by the projective reflections $(\sigma_s)_{s \in S}$ in the facets of $P$. Then the following hold:
\begin{enumerate}[label=(\roman*)]
\item the homomorphism $\sigma : W_P \rightarrow \SL^{\pm}(V)$ defined by $\sigma(s)=\sigma_s$ is an isomorphism onto $\Gamma_P$;
\item the group $\G_P$ is a discrete projective reflection group;
\item the union $\mathscr{C}_P$ of the $\Gamma_P$-translates of $P$ is a convex subset of $\PS(V)$;
\item if $\Omega_P$ is the interior of $\mathscr{C}_P$, then $\Gamma_P$ acts properly discontinuously on $\Omega_P$.
\end{enumerate}
\end{theorem}

\subsection{From Cartan matrices to mirror polytopes}\label{subsec:CartanToPolytope}

Given a Cartan matrix $\Cart_S$, there is a simple process to build a canonical mirror polytope $\Delta_\Cart$ such that $\Cart_{\Delta_\Cart} = \Cart_S$. In this construction, $\Delta_\Cart$ will be a simplex of dimension $\# S-1$.

\medskip

Let $V=\R^S$. We denote by $(e_s)_{s \in S}$ the canonical basis of $V$ and $(e^*_s)_{s \in S}$ its dual basis. We set $\alpha_s := e^*_s$ and $v_s := \sum_t \Cart_{t,s} \, e_t$, i.e., $v_s$ is the $s$-column vector of $\Cart_S$. Hence, by taking $\Delta$ to be (the projectivization of) the negative quadrant in $\PS(V)$ and $\sigma_s = \mathrm{Id} - \alpha_s \otimes v_s$ to be the reflection in the facet $\Delta \cap \PS( \mathrm{Ker}\,\alpha_s )$, we obtain a mirror polytope $\Delta_\Cart$ whose underlying polytope is a simplex of dimension $\# S-1$. We call $\Delta_\Cart$ the \emph{mirror simplex associated with $\Cart_S$}.

\medskip

In the case where $W_S$ is any Coxeter group and $\Cart_S = \mathsf{C}_W$, the mirror simplex associated with $\Cart_S$ is called the \emph{Tits simplex associated with $W_S$}\index{Tits simplex} and denoted by $\Delta_W$. The corresponding representation $\sigma : W_S \to \PS(\R^S)$ is dual to Tits geometric representation described in \cite{Bourbaki_group_456}. 
\medskip

For example, if $W$ is spherical (resp. irreducible affine), then $\Delta_W$ gives rise to the classical tiling of $\S (V)$ (resp. of an affine chart) with $\G_{\Delta_W}$ in the isometry group of the sphere (resp. Euclidean space). If $W$ is Lannér (resp. quasi-Lannér), then $\Omega_{\Delta_W}$ is the projective model of the hyperbolic space and $\Delta_W$ (resp. $\Delta_W \cap \Omega_{\Delta_W}$) is the hyperbolic Coxeter polytope whose Coxeter group is $W$.

\medskip
 
By the Perron--Frobenius theorem, an irreducible Cartan matrix $\Cart$ has a simple eigenvalue $\lambda_\Cart$ which corresponds to an eigenvector with positive entries and has the smallest modulus among the eigenvalues of $\Cart$. We say that $\Cart$ is of \emph{positive}\index{Cartan matrix!positive type}, \emph{zero}\index{Cartan matrix!zero type} or \emph{negative type}\index{Cartan matrix!negative type} when $\lambda_\Cart$ is positive, zero or negative, respectively. For example, the Gram matrix of a hyperbolic polytope of finite volume is always of negative type. Now, the following is a generalization of Theorem \ref{thm:gram} to the projective setting.

\begin{theorem}[Vinberg {\cite[Cor.~1]{bible}}]\label{thm:Cartan}
Let $\Cart$ be a Cartan matrix of size $N \times N$. Assume that $\Cart$ is irreducible, of negative type and of rank $d+1$. Then there exists a unique mirror $d$-polytope $P$, up to automorphism of $\S(\R^{d+1})$, such that $\Cart_{P} = \Cart$.
\end{theorem}

\begin{remark}
Theorem \ref{thm:Cartan} is not explicitly stated in \cite[Cor.~1]{bible} for non-Coxeter polytopes, but may be proved from \cite[Prop.~13 \& 15]{bible}. 
\end{remark}

\subsection{Anosov reflection groups}

Anosov representations are discrete representations of word-hyperbolic groups into semisimple Lie group with good dynamical properties. They have received a lot of attention and have been much studied recently (see e.g. \cite{labourie_anosov,Guichard_Wienhard_Anosov} for the definition of Anosov representation). But examples of Anosov representations of word hyperbolic groups, which are more complicated than free groups and surface groups, into Lie group of higher rank are less known. The following theorem, which generalizes Theorem \ref{thm:DH}, tells us that any infinite, word hyperbolic, irreducible Coxeter group admits Anosov representations. 

\begin{theorem}[{\cite[Cor.~1.18]{5A_cg_cc}}]
Let $P$ be a projective Coxeter polytope of $\PS(V)$ with Coxeter group $W_S$. Suppose that $W_S$ is word-hyperbolic. Then the following are equivalent:
\begin{itemize}
\item the representation $\sigma : W_S \rightarrow \SL^{\pm}(V)$ defined by $\sigma(s)=\sigma_s$ is $P_1$-Anosov (i.e., Anosov with respect to the stabilizer of a line in $V$);
\item $\Cart_{s,t} \Cart_{t,s} > 4$ for all $s\neq t$ with $m_{s,t} = \infty$.
\end{itemize}
\end{theorem}

\begin{remark}
Anosov reflection groups in $\mathrm{O}(p,q)$ can be used to give a new proof of Theorem \ref{thm:moussong} (Moussong's hyperbolicity criterion); see \cite{anosov_on_hpq,LM_quasi_fuchs}.
\end{remark}

\begin{theorem}[{\cite[Th.~A]{LM_quasi_fuchs}}]
In dimension $d=4, \ldots, 8$, there exists a projective Coxeter polytope of $\PS(\R^{d+2})$ with Coxeter group $W_S$ such that:
\begin{itemize}
	\item the group $W_S$ is word-hyperbolic and its boundary is a $(d-1)$-sphere;
	\item the image of the representation $\sigma : W_S \rightarrow \SL^{\pm}(\R^{d+2})$ defined by $\sigma(s)=\sigma_s$ lies in $\mathrm{O}_{d,2} (\R)$;
	\item the representation $\sigma : W_S \rightarrow \mathrm{O}_{d,2} (\R)$ is $P$-Anosov, where $P$ is the stabilizer of an isotropic line;
	\item the group $W_S$ is not quasi-isometric to $\Hb^d$.
\end{itemize}
\end{theorem}


\subsection{Convex cocompact projective reflection groups}

An infinite discrete subgroup $\Gamma$ of $\SL^{\pm}(V)$ is \emph{convex cocompact in $\PS(V)$}\index{convex cocompact in $\PS(V)$} if it acts properly discontinuously on some properly convex open subset $\Omega$ of $\PS(V)$ and cocompactly on a nonempty $\Gamma$-invariant closed convex subset $\mathscr{C}$ of $\Omega$ whose closure in $\PS(V)$ contains all accumulation points of all possible $\Gamma$-orbits $\Gamma \cdot y $ with $y \in \Omega$.

\medskip

The notion of Convex cocompactness in $\PS(V)$ introduced in \cite{anosov_on_RPV}, in some sense, generalizes that of Anosov representation, but it does not require that the group $\Gamma$ is word hyperbolic. There is also a simple characterization of convex cocompactness for projective reflection groups:

\begin{theorem}[{\cite[Th.~1.3]{5A_cg_cc}}]\label{thm:refl_gp_cc}
Let $P$ be a projective Coxeter polytope of $\PS(V)$ with infinite irreducible Coxeter group $W_S$, and $\sigma : W_S \rightarrow \SL^{\pm}(V)$ the representation defined by $\sigma(s)=\sigma_s$. If $\sigma(W_S)$ is convex cocompact in $\PS(V)$, then $W_S$ satisfies the following two conditions:
\begin{enumerate}[label=(\roman*)]
\item $S$ does not contain two orthogonal non-spherical subsets;
\item if $S$ contains an irreducible affine subset $T$ of rank $\geqslant 3$, then $W_T$ is of type $\widetilde{A}_k$ where $k = \# T -1$.
\end{enumerate} 
\end{theorem}

\begin{theorem}[{\cite[Th.~1.8]{5A_cg_cc}}]\label{thm:refl_gp_cc2}
Let $P$ be a projective Coxeter polytope of $\PS(V)$ with infinite irreducible Coxeter group $W_S$, $\Cart_S = (\Cart_{s,t})_{s,t \in S}$ the Cartan matrix of $P$, and $\sigma : W_S \rightarrow \SL^{\pm}(V)$ the representation defined by $\sigma(s)=\sigma_s$. If $W_S$ satisfies the conditions $(i)$ and $(ii)$ of Theorem \ref{thm:refl_gp_cc}, then the following are equivalent:
\begin{itemize}
\item $\sigma(W)$ is convex cocompact in $\PS(V)$;
\item for any irreducible standard subgroup $W_{T}$ of $W_S$ with $\varnothing \neq T \subset S$, the Cartan submatrix $\Cart_{T} := (\Cart_{s,t})_{s,t \in T}$ is not of zero type;
\item  $\det(\Cart_{T})\neq 0$ for all $T \subset S$ with $W_T$ of type~$\widetilde{A}_k$, $k\geqslant 1$.
\end{itemize}
\end{theorem}

As a result, any infinite, irreducible Coxeter group $W_S$ satisfying the conditions $(i)$ and $(ii)$ of Theorem \ref{thm:refl_gp_cc} admits projective reflection groups, which are convex cocompact in $\PS(\R^{N})$ with $N = \#S$ (see \cite[Th.~1.3]{5A_cg_cc}).

\subsection{Divisible and quasi-divisible domains}

Every properly convex open subset $\Omega$ of $\PS(V)$ admits a Hilbert metric $d_{\Omega}$ on $\Omega$ so that the group $\mathrm{Aut}(\Omega)$ of automorphisms of $\PS(V)$ preserving $\Omega$ acts on $\Omega$ by isometries for $d_{\Omega}$. A properly convex domain $\O$ is \emph{divisible}\index{divisible} (resp. \emph{quasi-divisible})\index{quasi-divisible} \emph{by $\Gamma$} if there exists a discrete subgroup $\G$ of $\Aut(\O)$ such that ${\O}/{\G}$ is compact (resp. of finite volume with respect to the Hausdorff measure induced by $d_{\Omega}$). (see e.g. \cite{marquis_handbook_hil_geo} for more details for the Hilbert metric and the Hausdorff measure). 

\medskip

In general, it is difficult to construct divisible or quasi-divisible domains with various properties. But, in small dimension, one can use \emph{perfect} or \emph{quasi-perfect} projective Coxeter polytopes to build such domains. We first introduce the definition of 2-perfect polytopes, which is slightly more general than that of perfect or quasi-perfect polytopes, and in Section \ref{sec:examples} we give some interesting examples of divisible domains.

\medskip

Let $P$ be a projective Coxeter polytope of $\PS(V)$ and $S$ the set of facets of $P$. Given a vertex $v$ of $P$, we denote by $S_v$ the set of facets that contain $v$. For any $s \in S_v$, the projective reflection $\sigma_s$ induces a projective reflection $\overline{\sigma}_s$ of the projective space $\PS( {V}/{\langle v\rangle} )$, where $\langle v\rangle$ is the subspace spanned by $v$ and ${V}/{\langle v\rangle}$ is the quotient vector space. The projection of $P$ to $\PS\left({V}/{\langle v\rangle}\right)$ with the reflections $(\overline{\sigma}_s)_{s \in S_v}$ define a projective Coxeter polytope $P_v$ of $\PS(V/{\langle v\rangle})$, called the \emph{link of $P$ at $v$}\index{projective Coxeter polytope!link}.

\begin{definition}
A projective Coxeter $d$-polytope $P$ is \emph{elliptic}\index{projective Coxeter polytope!elliptic} (resp. \emph{parabolic}\index{projective Coxeter polytope!parabolic}, resp. \emph{loxodromic})\index{projective Coxeter polytope!loxodromic} when each component of $\Cart_P$ is of positive type (resp. zero type, resp. negative type) and the rank of $\Cart_P$ is $d+1$ (resp. $d$, resp. $d+1$).
\end{definition}

\begin{remark}
If $P$ is elliptic, then $W_P$ is a spherical Coxeter group and $P$ is the Tits simplex associated with $W_P$. If $P$ is parabolic, then $W_P$ is an affine Coxeter group, $P$ is the Cartesian product of the Tits simplices associated with the irreducible components of $W_P$, and $\O_P$ is an affine chart of $\PS(V)$.
\end{remark}

A projective Coxeter polytope $P$ is \emph{perfect}\index{projective Coxeter polytope!perfect} (resp. \emph{quasi-perfect}\index{projective Coxeter polytope!quasi-perfect}, resp. \emph{2-perfect}\index{projective Coxeter polytope!2-perfect}) when all its vertex links are elliptic (resp. elliptic or parabolic, resp. perfect). For example, quasi-perfect Coxeter polytopes should be $2$-perfect.

\begin{remark}
By \cite[Prop.~26]{bible}, a perfect Coxeter polytope is either elliptic, parabolic or irreducible loxodromic.  
\end{remark}

Let $P$ be an irreducible loxodromic Coxeter polytope and $\Gamma_P$ the projective reflection group of $P$. Then $\Omega_P$ is a properly convex domain, hence it admits a Hilbert metric $d_{\Omega_P}$. By \cite[Th.~2]{bible}, a projective Coxeter polytope $P$ is perfect if and only if the action of $\G_P$ on $\O_P$ is cocompact. So, in this case, the domain $\O_P$ is divisible by $\G_P$.

\medskip

The action of $\G_P$ on $\O_P$ is said to be of \emph{finite covolume}\index{finite covolume} if $P \cap \Omega_P$ has finite volume with respect the Hausdorff measure $\mu_{\Omega_P}$ induced by $d_{\Omega_P}$, and \emph{geometrically finite}\index{geometrically finite} if $\mu_{\Omega_P} (P \cap \C(\Lambda_{\Omega_P})) < \infty$, where $\Lambda_P$ is the limit set of $\Gamma_P$ and $\C(\Lambda_P)$ is the convex hull of $\Lambda_P$ of $\Omega_P$ (see \cite{cox_in_hil} for more details).

\begin{theorem}[{\cite[Th.~A]{cox_in_hil}}]
Let $P$ be an irreducible, loxodromic, 2-perfect Coxeter polytope of $\PS(V)$. Then the action of $\G_P$ on $\O_P$ is always geometrically finite, and 
\begin{itemize}
  \item $\G_P$ is of finite covolume if and only if $P$ is quasi-perfect;
  
  \item $\G_P$ is convex cocompact in $\PS(V)$ if and only if all the vertex links of $P$ are elliptic or loxodromic.
\end{itemize}
\end{theorem}

\subsection{Cocompact action of Coxeter groups}

There are many examples of discrete Coxeter subgroups of $\mathrm{SL}^{\pm}(V)$ other than projective reflection groups. However, if a Coxeter group $\Gamma$ divides a properly convex domain, then $\Gamma = \Gamma_P$ for some projective Coxeter polytope $P$:

\begin{theorem}[Davis {\cite[Prop.~10.9.7]{davis_book} and Charney--Davis \cite{charney_davis_when}; see \cite[Lem.~5.4]{LM_quasi_fuchs}}]\label{thm:only_reflection}
Let $W$ be a Coxeter group, and let $\rho : W \rightarrow \mathrm{SL}^{\pm}(V)$ be a faithful representation. Suppose that there exists a convex domain $\Omega$ divisible by $\rho(W)$. Then the following hold:
\begin{enumerate}[label=(\roman*)]
\item for each $s \in S$, the image $\rho(s)$ of $s$ is a projective reflection of $\mathbb{S}(V)$;
\item $\rho(W)$ is a projective reflection group generated by $(\rho(s))_{s\in S}$.
\end{enumerate}
\end{theorem}

\begin{remark}
It is an open question whether Theorem \ref{thm:only_reflection} still holds when the word "divisible" is replaced by "quasi-divisible".
\end{remark}

\section{Examples of projective reflection groups}\label{sec:examples}

The construction of projective reflection groups had led to several existence theorems in convex projective geometry. 

\medskip

A properly convex domain $\Omega$ of $\PS(V)$ is \emph{decomposable}\index{decomposable} if a cone of $V$ lifting $\Omega$ is a non-trivial direct product of two smaller cones, and \emph{homogeneous}\index{homogeneous} if the group $\mathrm{Aut}(\Omega)$ acts transitively on $\Omega$. Since the homogeneous quasi-divisible domains are well-understood by \cite{homogeneous_vin_1,koecher}, only inhomogeneous ones are of interest to us. So, all properly convex domains in this section are assumed to be inhomogeneous and indecomposable.

\subsection{Kac--Vinberg's example}

\begin{figure}[ht!]
\centering
 \includegraphics[scale=.13]{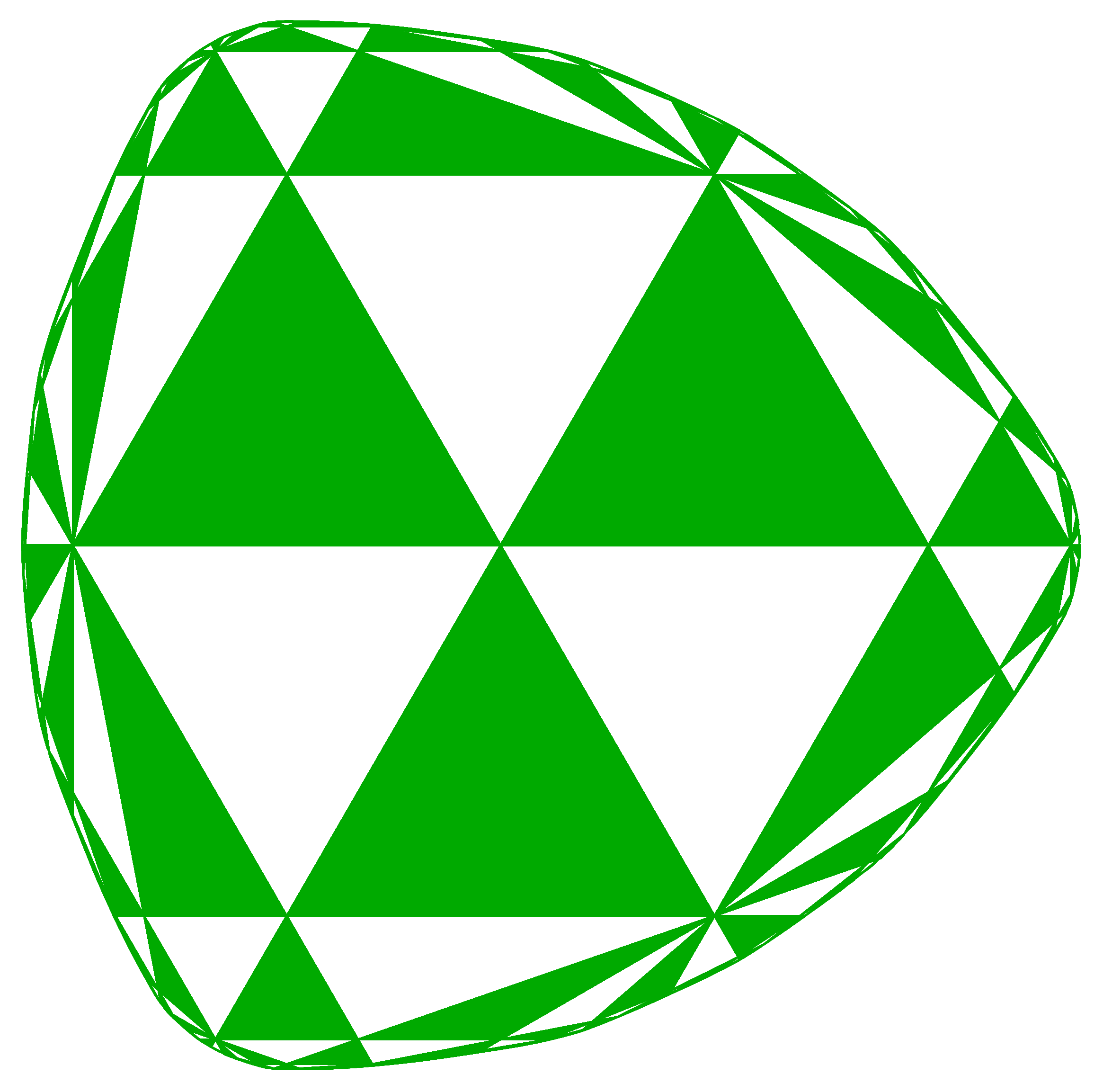} \quad\quad
 \includegraphics[scale=.13]{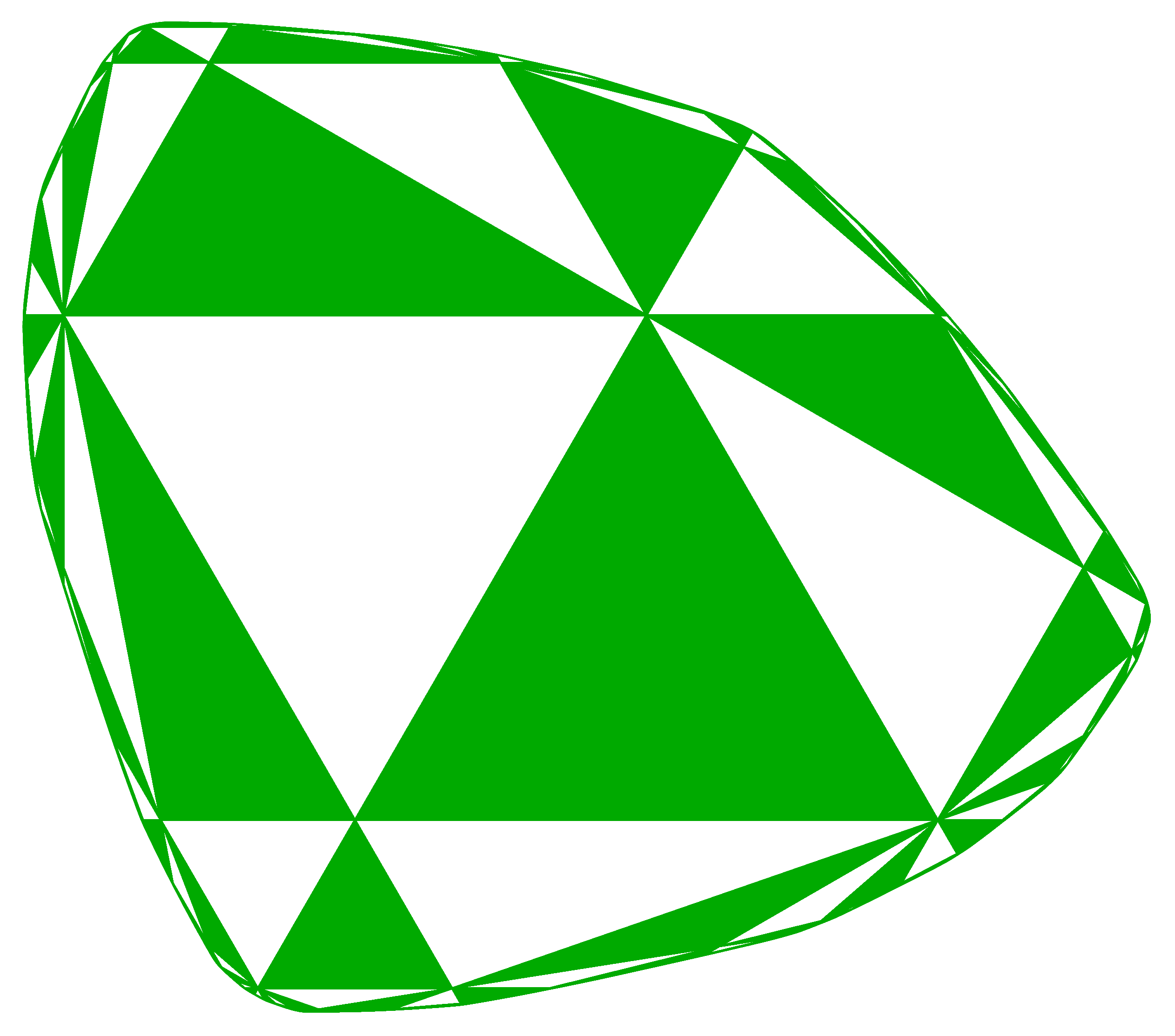} \quad\,
 \includegraphics[scale=.17]{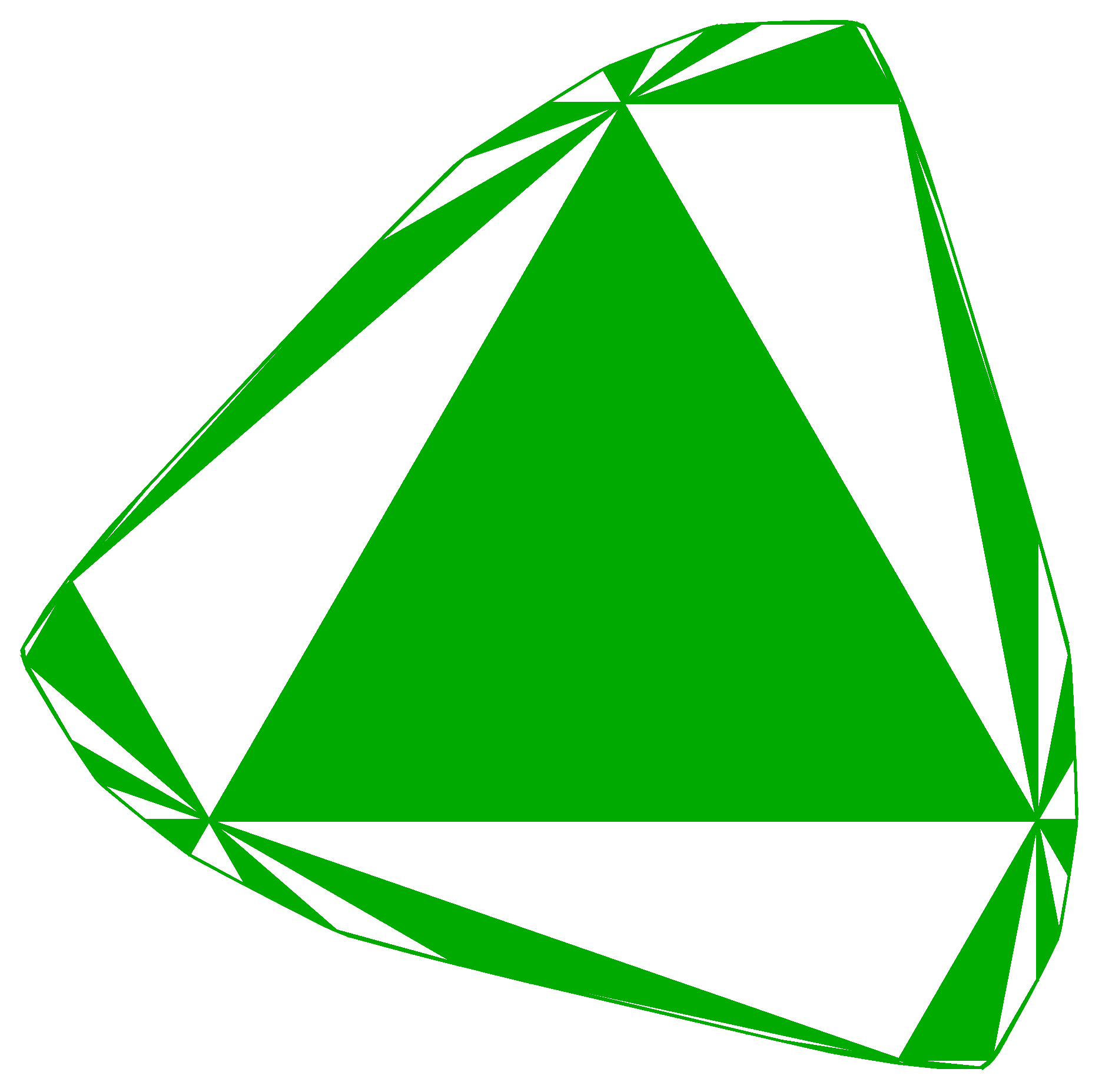}
 \caption{{Triangles with dihedral angles $\nicefrac{\pi}{3}$, $\nicefrac{\pi}{3}$ and $\nicefrac{\pi}{6}$ on the left, with dihedral angles $\nicefrac{\pi}{3}$, $\nicefrac{\pi}{4}$ and $\nicefrac{\pi}{6}$ on the center, and with dihedral angles $\nicefrac{\pi}{6}$, $\nicefrac{\pi}{6}$ and $\nicefrac{\pi}{6}$ on the right.}}\label{fig:Kac-Vinberg}
\end{figure}

The first example of divisible 2-domain which is not a hyperbolic plane was found by Kac and Vinberg \cite{kac_vinberg}. They used perfect projective Coxeter triangles $P$ with Cartan matrix $\Cart_P = (\Cart_{i,j})_{i,j=1,2,3}$ such that $(i)$ each entry $\Cart_{i,j}$ is an integer, $(ii)$ $\det(\Cart_P) < 0$, and $(iii)$ $\Cart_{1,2}\Cart_{2,3}\Cart_{3,1} \neq \Cart_{1,3}\Cart_{3,2}\Cart_{2,1}$. Here, the condition $(i)$ implies that the projective reflection group $\Gamma_P$ of $P$ is a subgroup of $\mathrm{SL}(3,\mathbb{Z})$, $(ii)$ implies that $\Cart_P$ is of negative type, and finally $(iii)$ implies that $\Omega_P$ is not a hyperbolic plane (see Figure \ref{fig:Kac-Vinberg}).

\begin{remark}
Let $\hat{\Gamma}_P$ be any finite-index torsion-free subgroup of $\Gamma_P$. Then $\hat{\Gamma}_P$ is an infinite index subgroup of $\mathrm{SL}(3,\mathbb{Z})$ and is Zariski dense in $\mathrm{SL}(3,\mathbb{R})$. In other words, $\hat{\Gamma}_P$ is a \emph{thin} surface group (see \cite{WhatIsThinGroup} for an introduction to thin groups).
\end{remark}

\subsection{Benoist's examples and more}

\begin{figure}[ht!]
\centering
 \includegraphics[scale=.25]{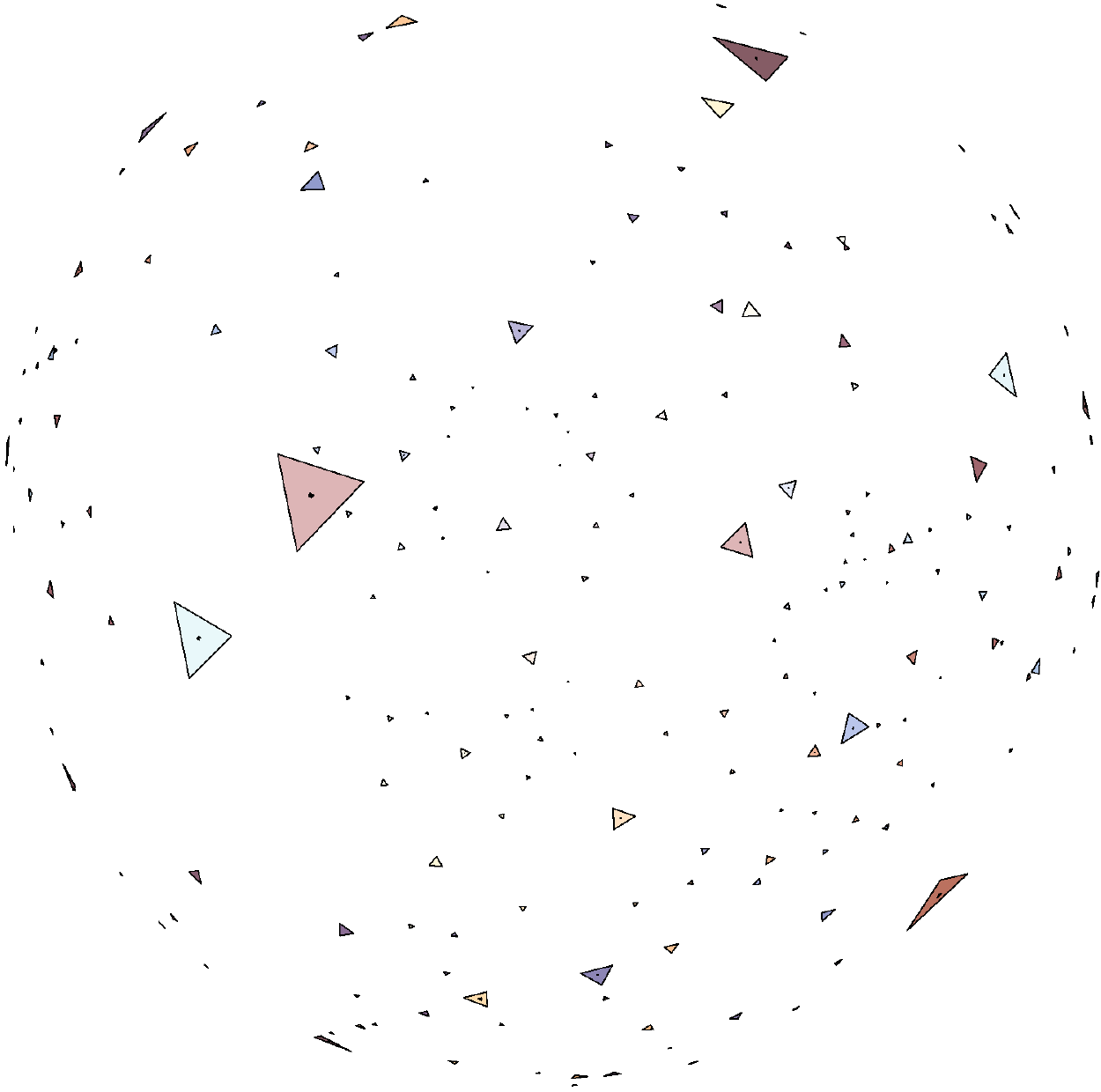} \quad\quad
 \includegraphics[scale=.25]{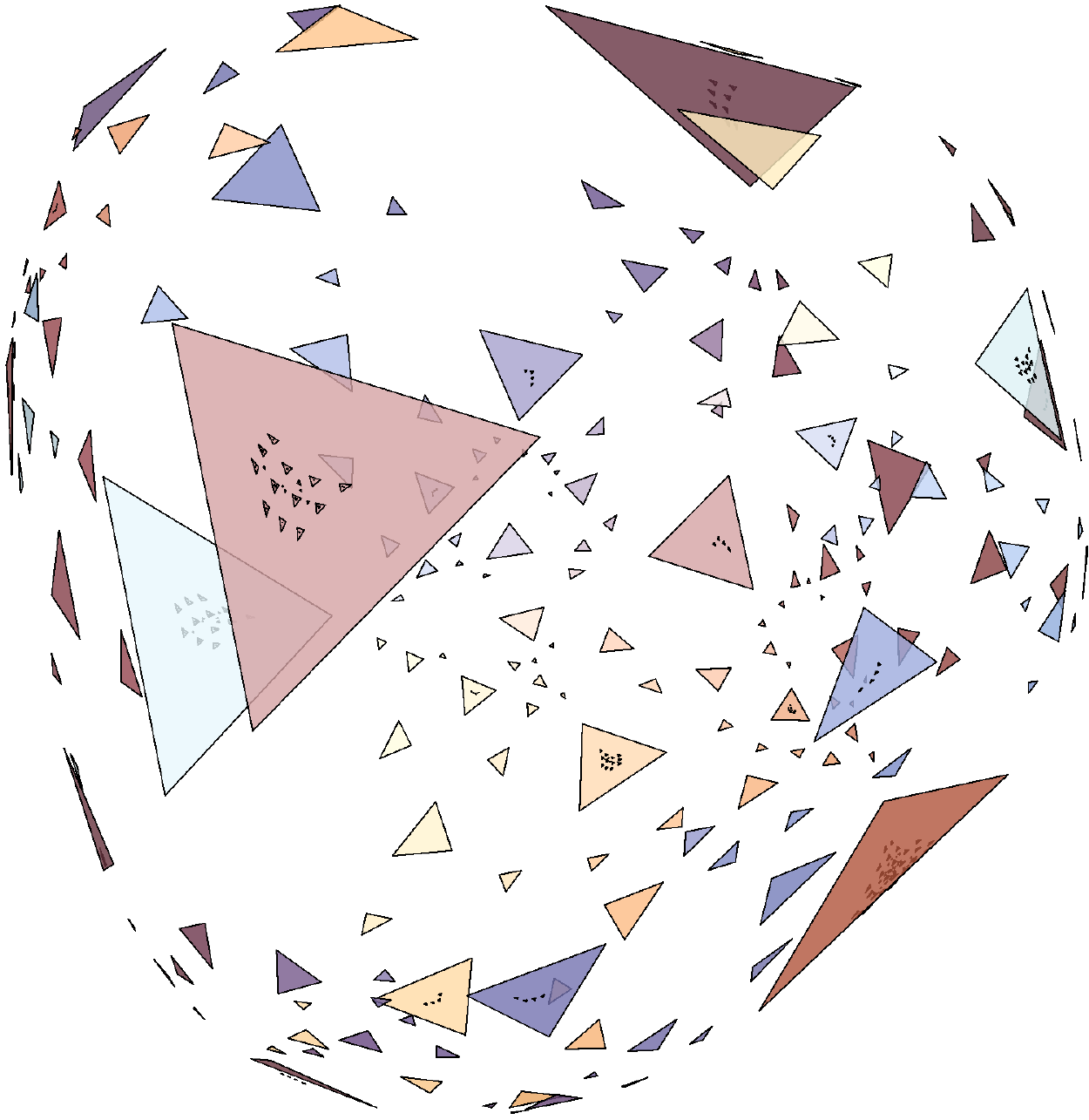} \quad\quad
 \includegraphics[scale=.25]{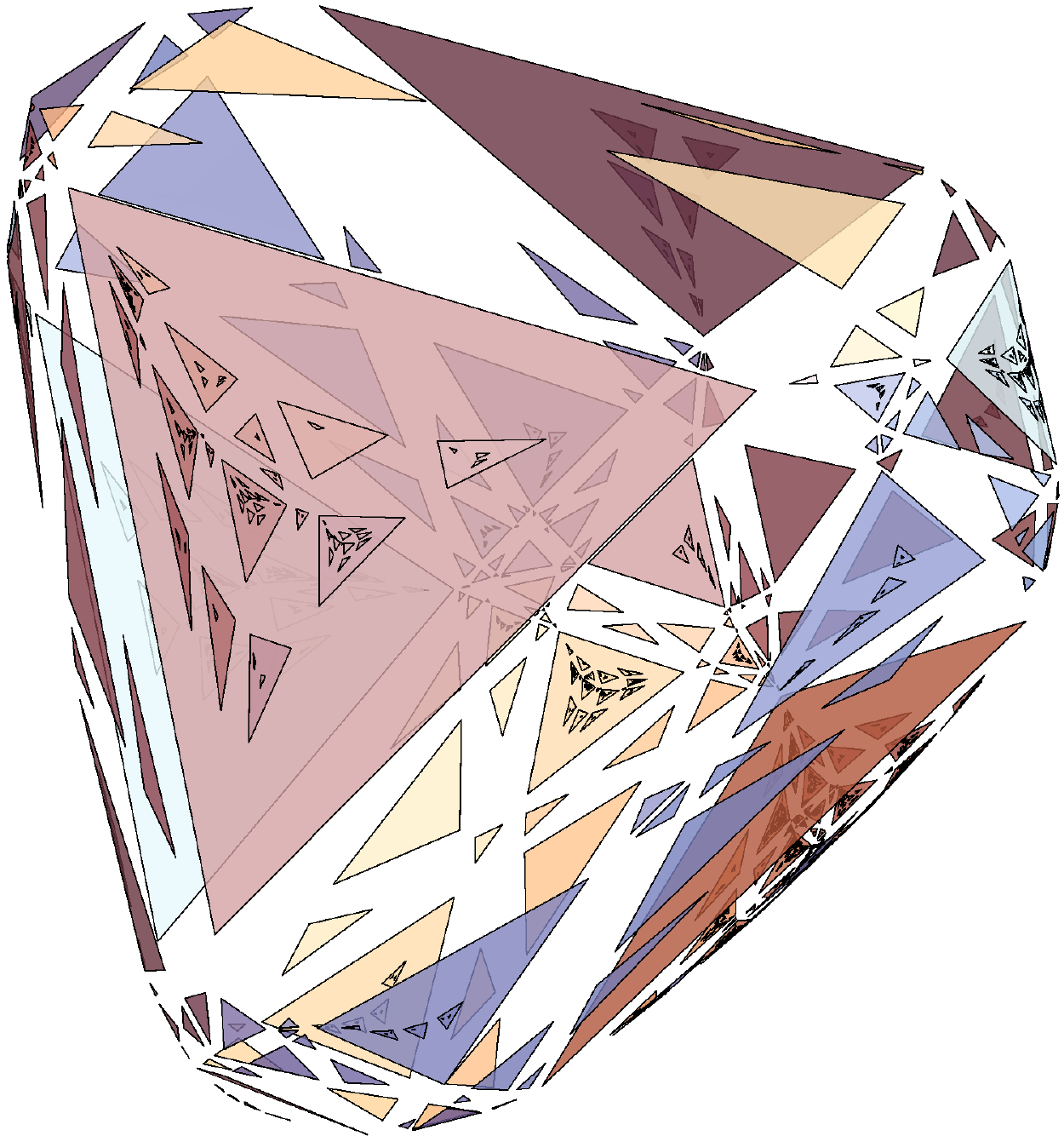}
 \caption{A collection of properly embedded triangles in the non-strictly convex divisible $3$-domains is colored. Each triangle is preserved by a subgroup of $\Gamma$, which is virtually $\mathbb{Z}^2$.}
 \label{fig:Benoist}
\end{figure}

The first known examples of divisible $d$-domains $\Omega$ which are not strictly convex were introduced by Benoist \cite{cd4} in dimension $d = 3, \dotsc, 7$ (see Figure \ref{fig:Benoist}). In such examples, the discrete group $\Gamma$ which divides $\Omega$ is relatively hyperbolic with respect to virtual $\Z^{d-1}$. Later, different examples of non-strictly convex divisible $d$-domains were found in \cite{CLM_dehn_fill} in dimension $d = 4, \dotsc, 8$, and the group $\Gamma$ dividing such $d$-domain is relatively hyperbolic with respect to a collection of virtually free abelian subgroup of rank $< {d-1}$.  Except in dimension 3 (see \cite{BDL_3d_geometrization}), all the known examples were built from projective reflection groups.  

\begin{remark}
A generalization of Thurston's hyperbolic Dehn filling theorem to the projective setting led to the examples in \cite{CLM_dehn_fill}.
\end{remark}

In \cite{benoist_qi}, Benoist found the first example of word-hyperbolic group $\Gamma$, not quasi-isometric to the hyperbolic space, that divides a properly convex 4-domain $\Omega$, again using projective reflection groups. Since $\Gamma$ is word hyperbolic, $\Omega$ should be strictly convex by \cite{convexe_div_1}. Shortly after, Kapovich \cite{kapo_qi} found examples in any dimension $d\geqslant 4$, using Gromov--Thurston manifolds \cite{gromov_thurston_manifold}.

\section{Hitchin component of polygon groups}

Let $P$ be a compact hyperbolic polygon with dihedral angles $\nicefrac{\pi}{m_1}, \dotsc, \nicefrac{\pi}{m_k}$, and let $W$ be the Coxeter group of $P$. The conjugacy classes of discrete and faithful representations of $W$ to $\PGL(2,\R)$ form a connected component $\mathcal{T}$ of $\chi(W,\PGL(2,\R)):= \mathrm{Hom}(W,\PGL(2,\R)) / \PGL(2,\R)$, i.e., the space of conjugacy classes of representations of $W$ to $\PGL(2,\R)$.

\medskip

For any $n \geqslant 2$, there is a unique irreducible representation $\kappa : \PGL(2,\R) \rightarrow \PGL(n,\R)$ up to conjugation. This gives rise to an embedding:
$$ \mathcal{T} \rightarrow \chi(W,\PGL(n,\R)) := \mathrm{Hom}(W,\PGL(n,\R)) / \PGL(n,\R)$$
The image of this embedding is called the \emph{Fuchsian locus}\index{Fuchsian locus} and the component of $\chi(W,\PGL(n,\R))$ containing the Fuchsian locus is the \emph{Hitchin component}\index{Hitchin component} $\Hit(W, \PGL(n,\R))$. A representation $\rho : W \rightarrow \PGL(n,\R)$ is called a \emph{Hitchin representation}\index{Hitchin representation} if its $\PGL(n,\R)$-conjugacy class is an element of $\Hit(W, \PGL(n,\R))$.

\begin{theorem}[{\cite[Th.~1.1 \& 1.2]{ALS}}]\label{thm:Hithcin}
Let $P$ be a compact hyperbolic polygon with dihedral angles $\nicefrac{\pi}{m_1}, \dotsc, \nicefrac{\pi}{m_k}$, and $W$ its Coxeter group. Then each Hitchin representation in $\Hit(W,\PGL(n,\R))$ is discrete and faithfull, and $\Hit(W,\PGL(n,\R))$ is an open cell of dimension 
$$ -(n^2-1) + \sum_{\ell=2}^n \sum_{i=1}^k \left \lfloor{ \ell \left( 1-\frac{1}{m_i} \right)}\right \rfloor,$$
where $\lfloor x \rfloor$ denotes the biggest integer not bigger than $x$.
\end{theorem}

For example, the $\PGL(2m,\R)$ (resp. $\PGL(2m+1,\R)$) Hitchin component of the Coxeter group associated with a right-angled hyperbolic $k$-gon ($k \geqslant 5$) is an open cell of dimension $(k-4)m^2+1$ (resp. $(k-4)(m^2+m)$).

\begin{remark}
In the case of $n=2$ (resp. $n=3$), Theorem \ref{thm:Hithcin} was proved by Thurston \cite{CoursdeThurston} (resp. Choi--Goldman \cite{ChoiGoldman2005}). 
\end{remark}

\begin{remark}
Let $\rho$ be any Hitchin representation in $\Hit(W,\PGL(n,\R))$. In the case $n \geqslant 4$, the image of each generator of $W$ should be an involution but not a projective reflection, hence $\rho$ is not a projective reflection group.
\end{remark}

\section{Properly discontinuous affine groups}

\subsection{Auslander's conjecture and Milnor's question}

In the 1960s, Auslander raised the following conjecture: 

\begin{conj}[Auslander \cite{Auslander1964}]
Every discrete subgroup $\Gamma$ of the affine group $\mathrm{Aff}(\R^d)$ which acts properly discontinuously and cocompactly on $\R^d$ is virtually solvable.
\end{conj}

In the 1970s, Milnor asked if Auslander's conjecture still holds without the condition that the action is cocompact:

\begin{ques}[Milnor \cite{Milnor1977}]
Is every discrete subgroup $\Gamma$ of $\mathrm{Aff}(\R^d)$ which acts properly discontinuously on $\R^d$ virtually solvable?
\end{ques}

In 1983, Fried and Goldman \cite{Fried_Goldman} showed that Auslander's conjecture is true in dimension $3$, and Margulis answered Milnor's question negatively:
\begin{theorem}[Margulis \cite{Margulis_affine_1,Margulis_affine_2}]
There exists a properly discontinuous affine action of the free group on two generators on $\R^3$.
\end{theorem}

Even if some progress have been made over the years towards Auslander's conjecture (see e.g. \cite{Tomanov_2016,AMS_2012_RF} for a proof assuming $d \leqslant 6$ and \cite{Goldman_Kamishima,Tomanov1990,AMS_2010} for a proof assuming the linear part is contained in a particular class of semisimple Lie subgroups), Auslander's conjecture is still open.

\medskip

Back to Milnor's question, the existence and property of properly discontinuous affine action of free groups on $\R^n$ have been actively studied (see e.g. \cite{Drumm1992,CDG_2016,DGK_arc_complex,GLM_affine,AMS_jdg,Smilga_1,Smilga_3} or the survey \cite{DDGS_survey}).

\subsection{Properly discontinuous affine Coxeter groups}

Before the following theorem, properly discontinuous affine actions by non-virtually solvable non-free groups were unknown.

\begin{theorem}[Danciger--Guéritaud--Kassel {\cite[Th.~1.1]{DGK_affine_Cox}}]
Any right-angled Coxeter group of rank $k$ admits a properly discontinuous affine action on $\R^{k(k-1)/2}$.
\end{theorem}

\begin{remark}
The action preserves a bilinear form and in some particular cases, one can find much smaller affine space on which the Coxeter groups acts (see \cite[Prop.~1.6]{DGK_affine_Cox}).
\end{remark}


\newcommand{\etalchar}[1]{$^{#1}$}

\printindex

\end{document}